\documentclass[12pt]{amsart}
\usepackage{amssymb}
\usepackage[all]{xy}
\usepackage[toc,page,title,titletoc,header]{appendix}
\usepackage{xcolor}

\begin{document}
\theoremstyle{plain}
\newtheorem{thm}{Theorem}[section]
\newtheorem*{thm1}{Theorem 1}
\newtheorem*{thm1.1}{Theorem 1.1}
\newtheorem*{thmM}{Main Theorem}
\newtheorem*{thmA}{Theorem A}
\newtheorem*{thm2}{Theorem 2}
\newtheorem{lemma}[thm]{Lemma}
\newtheorem{lem}[thm]{Lemma}
\newtheorem{cor}[thm]{Corollary}
\newtheorem{pro}[thm]{Proposition}
\newtheorem{propose}[thm]{Proposition}
\newtheorem{variant}[thm]{Variant}
\theoremstyle{definition}
\newtheorem{notations}[thm]{Notations}
\newtheorem{rem}[thm]{Remark}
\newtheorem{rmk}[thm]{Remark}
\newtheorem{rmks}[thm]{Remarks}
\newtheorem{defi}[thm]{Definition}
\newtheorem{exe}[thm]{Example}
\newtheorem{claim}[thm]{Claim}
\newtheorem{ass}[thm]{Assumption}
\newtheorem{prodefi}[thm]{Proposition-Definition}
\newtheorem{que}[thm]{Question}
\newtheorem{con}[thm]{Conjecture}
\newtheorem*{assa}{Assumption A}

\newtheorem*{dmlcon}{Dynamical Mordell-Lang Conjecture}
\newtheorem*{condml}{Dynamical Mordell-Lang Conjecture}
\numberwithin{equation}{section}
\newcounter{elno}                
\def\points{\list
{\hss\llap{\upshape{(\roman{elno})}}}{\usecounter{elno}}}
\let\endpoints=\endlist
\newcommand{\SH}{\rm SH}
\newcommand{\Tan}{\rm Tan}
\newcommand{\res}{\rm res}
\newcommand{\Om}{\Omega}
\newcommand{\om}{\omega}
\newcommand{\la}{\lambda}
\newcommand{\mc}{\mathcal}
\newcommand{\mb}{\mathbb}
\newcommand{\surj}{\twoheadrightarrow}
\newcommand{\inj}{\hookrightarrow}
\newcommand{\zar}{{\rm zar}}
\newcommand{\Exc}{\rm Exc}
\newcommand{\an}{{\rm an}}
\newcommand{\Con}{{\rm con}}
\newcommand{\red}{{\rm red}}
\newcommand{\codim}{{\rm codim}}
\newcommand{\Supp}{{\rm Supp}}
\newcommand{\rank}{{\rm rank}}
\newcommand{\Ker}{{\rm Ker \ }}
\newcommand{\Pic}{{\rm Pic}}
\newcommand{\Der}{{\rm Der}}
\newcommand{\Div}{{\rm Div}}
\newcommand{\Hom}{{\rm Hom}}
\newcommand{\im}{{\rm im}}
\newcommand{\Spec}{{\rm Spec \,}}
\newcommand{\Nef}{{\rm Nef \,}}
\newcommand{\Frac}{{\rm Frac \,}}
\newcommand{\Sing}{{\rm Sing}}
\newcommand{\sing}{{\rm sing}}
\newcommand{\reg}{{\rm reg}}
\newcommand{\Char}{{\rm char}}
\newcommand{\Tr}{{\rm Tr}}
\newcommand{\ord}{{\rm ord}}
\newcommand{\id}{{\rm id}}
\newcommand{\NE}{{\rm NE}}
\newcommand{\Gal}{{\rm Gal}}
\newcommand{\Min}{{\rm Min \ }}
\newcommand{\Max}{{\rm Max \ }}
\newcommand{\Alb}{{\rm Alb}\,}
\newcommand{\GL}{{\rm GL}\,}        
\newcommand{\PGL}{{\rm PGL}\,}
\newcommand{\Bir}{{\rm Bir}}
\newcommand{\Aut}{{\rm Aut}}
\newcommand{\End}{{\rm End}}
\newcommand{\Per}{{\rm Per}\,}
\newcommand{\ie}{{\it i.e.\/},\ }
\newcommand{\niso}{\not\cong}
\newcommand{\nin}{\not\in}
\newcommand{\soplus}[1]{\stackrel{#1}{\oplus}}
\newcommand{\by}[1]{\stackrel{#1}{\rightarrow}}
\newcommand{\longby}[1]{\stackrel{#1}{\longrightarrow}}
\newcommand{\vlongby}[1]{\stackrel{#1}{\mbox{\large{$\longrightarrow$}}}}
\newcommand{\ldownarrow}{\mbox{\Large{\Large{$\downarrow$}}}}
\newcommand{\lsearrow}{\mbox{\Large{$\searrow$}}}
\renewcommand{\d}{\stackrel{\mbox{\scriptsize{$\bullet$}}}{}}
\newcommand{\dlog}{{\rm dlog}\,}    
\newcommand{\longto}{\longrightarrow}
\newcommand{\vlongto}{\mbox{{\Large{$\longto$}}}}
\newcommand{\limdir}[1]{{\displaystyle{\mathop{\rm lim}_{\buildrel\longrightarrow\over{#1}}}}\,}
\newcommand{\liminv}[1]{{\displaystyle{\mathop{\rm lim}_{\buildrel\longleftarrow\over{#1}}}}\,}
\newcommand{\norm}[1]{\mbox{$\parallel{#1}\parallel$}}
\newcommand{\boxtensor}{{\Box\kern-9.03pt\raise1.42pt\hbox{$\times$}}}
\newcommand{\into}{\hookrightarrow}
\newcommand{\image}{{\rm image}\,}
\newcommand{\Lie}{{\rm Lie}\,}      
\newcommand{\CM}{\rm CM}
\newcommand{\sext}{\mbox{${\mathcal E}xt\,$}}  
\newcommand{\shom}{\mbox{${\mathcal H}om\,$}}  
\newcommand{\coker}{{\rm coker}\,}  
\newcommand{\sm}{{\rm sm}}
\newcommand{\pgcd}{\text{pgcd}}
\newcommand{\trd}{\text{tr.d.}}
\newcommand{\tensor}{\otimes}
\newcommand{\hotimes}{\widehat{\otimes}}

\renewcommand{\iff}{\mbox{ $\Longleftrightarrow$ }}
\newcommand{\supp}{{\rm supp}\,}
\newcommand{\ext}[1]{\stackrel{#1}{\wedge}}
\newcommand{\onto}{\mbox{$\,\>>>\hspace{-.5cm}\to\hspace{.15cm}$}}
\newcommand{\propsubset}
{\mbox{$\textstyle{
\subseteq_{\kern-5pt\raise-1pt\hbox{\mbox{\tiny{$/$}}}}}$}}
\newcommand{\sA}{{\mathcal A}}
\newcommand{\sB}{{\mathcal B}}
\newcommand{\sC}{{\mathcal C}}
\newcommand{\sD}{{\mathcal D}}
\newcommand{\sE}{{\mathcal E}}
\newcommand{\sF}{{\mathcal F}}
\newcommand{\sG}{{\mathcal G}}
\newcommand{\sH}{{\mathcal H}}
\newcommand{\sI}{{\mathcal I}}
\newcommand{\sJ}{{\mathcal J}}
\newcommand{\sK}{{\mathcal K}}
\newcommand{\sL}{{\mathcal L}}
\newcommand{\sM}{{\mathcal M}}
\newcommand{\sN}{{\mathcal N}}
\newcommand{\sO}{{\mathcal O}}
\newcommand{\sP}{{\mathcal P}}
\newcommand{\sQ}{{\mathcal Q}}
\newcommand{\sR}{{\mathcal R}}
\newcommand{\sS}{{\mathcal S}}
\newcommand{\sT}{{\mathcal T}}
\newcommand{\sU}{{\mathcal U}}
\newcommand{\sV}{{\mathcal V}}
\newcommand{\sW}{{\mathcal W}}
\newcommand{\sX}{{\mathcal X}}
\newcommand{\sY}{{\mathcal Y}}
\newcommand{\sZ}{{\mathcal Z}}
\newcommand{\A}{{\mathbb A}}
\newcommand{\B}{{\mathbb B}}
\newcommand{\C}{{\mathbb C}}
\newcommand{\D}{{\mathbb D}}
\newcommand{\E}{{\mathbb E}}
\newcommand{\F}{{\mathbb F}}
\newcommand{\G}{{\mathbb G}}
\newcommand{\HH}{{\mathbb H}}
\newcommand{\I}{{\mathbb I}}
\newcommand{\J}{{\mathbb J}}
\newcommand{\M}{{\mathbb M}}
\newcommand{\N}{{\mathbb N}}
\renewcommand{\P}{{\mathbb P}}
\newcommand{\Q}{{\mathbb Q}}
\newcommand{\R}{{\mathbb R}}
\newcommand{\T}{{\mathbb T}}
\newcommand{\U}{{\mathbb U}}
\newcommand{\V}{{\mathbb V}}
\newcommand{\W}{{\mathbb W}}
\newcommand{\X}{{\mathbb X}}
\newcommand{\Y}{{\mathbb Y}}
\newcommand{\Z}{{\mathbb Z}}
\newcommand{\bk}{{\mathbf{k}}}
\newcommand{\bt}{{\mathbf{t}}}

\newcommand{\fix}{\mathrm{Fix}}
\newcommand{\Fun}{\mathbf{Fun}}

\title[]{Continuous functions on Berkovich spaces}

\author{Junyi Xie}

\address{Beijing International Center for Mathematical Research, Peking University, Beijing 100871, China}

\email{xiejunyi@bicmr.pku.edu.cn}

\thanks{The author is supported by the NSFC Grant No.12271007 and 
partially supported by project ``Fatou'' ANR-17-CE40-0002-01}

\date{\today}

\bibliographystyle{plain}

\maketitle

\begin{abstract}Let $\bk$ be a perfect complete valued field with a nontrivial non-archimedean norm $|\cdot|$ and $\omega\in \bk$ with $0<|\omega|<1.$
Let $X$ be a reduced and normal $\bk$-analytic space. Then $\sO^{\circ}\simeq \lim\limits_{\leftarrow}\sO^{\circ}/\omega^n.$
\end{abstract}


\section{Introduction}
\subsection{A motivation}
Let $X$ be a complex manifold. Cauchy's convergence theorem says that a Cauchy sequence of continuous (resp. holomorphic) functions for the $L^{\infty}$-norm convergences to a continuous (resp. holomorphic) function on $X$.

Let us consider a sheaf version of this theorem. Denote by $\sC$ (resp. $\sO$) the sheaf of $\C$-valued continuous (resp. holomorphic) functions.
Set $\sF=\sC$ or $\sO.$ For every open subset $U\subseteq X$, denote by $\rho_U: \sF(U)\to [0,+\infty]$ the $L^{\infty}$-norm. Define a \emph{piecewise $\sF$-function on $U$} to be a collection $\{(U_i,F_i), i\in I\}$ where $U_i, i\in I$ is an open covering of $U$ and $F_i\in \sF(U_i).$ Denote by $\sP(U)$ the set of piecewise $\sF$-function on $U$.

For $F=\{(U_i,F_i), i\in I\},G=\{(V_j,G_j), j\in J\}\in \sP(U)$, define $$d(F,G):=\sup_{i\in I,j\in J}\rho((F_i|_{U_i\cap U_j}-G_j|_{U_i\cap U_j})).$$ The function $d$ is similar to a semi-metric, except that $d(F,F)$ may take positive value.
We may define the Cauchy sequences on $\sP(U)$ as usual, to be the sequences $F_n, n\geq 0$, for which, for every $\epsilon>0$, there is $N\geq 0$ such that for every $m,n\geq N$, $d(F_n,F_m)\leq \epsilon.$
We define the equivalence relation of the Cauchy sequences in the usual way.
Define $\hat{\sF}$ to be the presheaf on $X$ sending $U$ to the equivalence classes of Cauchy sequences in $\sP(U)$. Easy to see that $\sF$ is a sub-presheaf of $\hat{\sF}.$

When $\sF=\sC$, we have $\hat{\sF}=\hat{\sC}=\sC.$ This fact can be viewed as a sheaf version of Cauchy's convergence theorem for $\sC.$
But this sheaf version is not true for $\sO$ when $\dim X\geq 1.$ Observe that every bounded continuous function can be uniformly approximated by piecewise locally constant functions (hence by piecewise holomorphic functions). So $\hat{\sO}(X)$ contains all bounded continuous functions on $X$. In fact, it is easy to see that $\sC$ is exactly
the sheafification of $\hat{\sF}$. This gives us a way to define the sheaf of continuous functions from the sheaf of holomorphic functions.

In non-archimedean geometry (in the sence of Berkovich), we a priori have a structure sheaf which is the analogy of the sheaf of holomorphic functions. Applying the above construction, we get an analogy of the sheaf of continuous functions. One may ask what does it look like?  In particular, are there ``continuous functions" which are not analytic? In this paper we show that under mild conditions, the sheaf version of Cauchy's convergence theorem holds for the structure sheaf. In other words, all ``continuous functions" in this sense are analytic. We will also study the analytic and continuous functions on subsets of Berkovich spaces. In particular, there are continuous functions on some closed subsets which are not analytic. 

In the sequels to the papers, we will apply the results in this paper to study non-archimedean dynamics. In particular, Theorem \ref{thmisolim} will be used to 
generate \cite[Appendix B]{Xie2019} to the global setting.

\subsection{The statement}
Denote by $\bk$ a complete valued field with a nontrivial non-archimedean norm $|\cdot |.$ Denote by $\bk^{\circ}:=\{f\in \bk|\,\,|f|\leq 1\}$ the valuation ring and $\bk^{\circ\circ}:=\{f\in \bk|\,\,|f|< 1\}$ its maximal ideal. Denote by $\widetilde{\bk}:=\bk^{\circ}/\bk^{\circ\circ}$ the residue field. 

\medskip

Let $X$ be a reduced $\bk$-analytic space. We denote by $\rho_X$ the spectral norm on $X$. Some times, we write $\rho$ for $\rho_X$ for the simplicity.
Let $\sO$ (or $\sO_X$ when we want to emphasize the space $X$) be the structure sheaf on $X$ w.r.t the $G$-topology.
Let $\sO^{\circ}(X)$ the subsheaf of $\sO$ satisfying $\sO^{\circ}(U)=\sO(U)^{\circ}$ for every analytic subdomain $U$ of $X$. Let $\omega$ be any element in $\bk$ satisfying $|\omega|\in (0,1).$
Let $\sO^{\circ}/\omega$ be the cokernal of the morphism $\omega \times \cdot:\sO^{\circ}\to \sO^{\circ}.$
 We prove the following result.

\begin{thm}\label{thmisolim}Assume that $\bk$ if perfect and $X$ is normal. Let $\omega$ be any element in $\bk$ satisfying $|\omega|\in (0,1).$
Denote by $\alpha_n: \sO^{\circ}\to \sO^{\circ}/\omega^n$ the quotient morphism. Then the natural morphism $\alpha:= \lim\limits_{\leftarrow} \alpha_n: \sO^{\circ}\to \lim\limits_{\leftarrow}\sO^{\circ}/\omega^n$ is an isomorphism.
\end{thm}

As pointed out by the referee of the first version of the paper, 
the rigid space version of this theorem has a simplified  proof using \cite[Theorem 2]{Bartenwerfer1978}. 
The rigid space version means that $X$ itself is strict and we only consider strict covers.
We will sketch his proof in Section \ref{subsectionanotherproof}.

\begin{rem}
It is easy to see that the morphism $\sO^{\circ}\to \lim\limits_{\leftarrow}\sO^{\circ}/\omega^n$ is an isomorphism if and only if the natural morphism $\sO\to \lim\limits_{\leftarrow}\sO/\omega^n\sO^{\circ}$ is an isomorphism.
\end{rem}

\begin{rem}
As Ruochuan Liu pointed out, when $\bk$ is of mixed characteristic,
there is an alternative proof of Theorem \ref{thmisolim} using the almost vanishing of the higher cohomologies of $\sO^{+}$ on affinoid perfectoid spaces \cite[Theorem 6.3]{Scholze2012} and  the Ax-Sen-Tate for rigid analytic spaces \cite[Proposition 8.2.1]{Kedlaya}. He also pointed out that when $\bk$ is of mixed characteristic, Theorem \ref{thmisolim} is implied by \cite[Theorem 10.3]{Hansen2020}.
\end{rem}

\begin{rem}\label{remnormalness}If $X$ is not normal, $\alpha$ can be non-surjective.
An example is as follows: Let $X:=\sM(\bk\{T_1,T_2\}/(T_2^2-T_1^3)).$ For $n\geq 0$,  let $U_n:=\{x\in X|\,\, |T_1(x)|\leq |\omega^{2n}|\}$ and $V_n:=\{x\in X|\,\, |T_1(x)|\geq |\omega^{2n}|\}.$
We note that for $x\in V_n$, $T_2/T_1\in \sO(V_1)$ and $|T_2/T_1(x)|=|T_1(x)|^{1/2}.$ So we have $T_2/T_1\in \sO(V_n)^{\circ}.$ Set $v_n:=\alpha_n(T_2/T_1)\in \sO^{\circ}/\omega^n(V_n).$ 
Since for every $x\in U_n\cap V_n$, $|T_2/T_1(x)|=|T_1(x)|^{1/2}=|\omega^{n}|$, we have
$v_n|_{U_n\cap V_n}=0\in \sO^{\circ}/\omega^n(V_1\cap V_2).$ So there is a unique section $t_n\in \sO^{\circ}/\omega^n(X)$ such that $t_n|_{V_n}=v_n$ and $t_n|_{U_n}=0.$
It is clear that for $m\geq n$, the image of $t_m$ in $\sO^{\circ}/\omega^n(X)$ is $t_n.$ Set $T:=\lim\limits_{\leftarrow}t_n\in \lim\limits_{\leftarrow}\sO^{\circ}/\omega^n.$
Let $o$ be the $\bk$-point in $X$ defined by the maximal ideal $(T_1, T_2)$ in $\bk\{T_1,T_2\}/(T_2^2-T_1^3).$ 
We note that on $Y:=X\setminus \{o\}$, we have $T|_Y=\alpha(T_2/T_1)$. But $T_2/T_1$ does not extend to a regular function on $X$. So we have $T\not\in \alpha(\sO^{\circ}(X)).$ 
\end{rem}

The quotient morphism $\alpha_{\omega}:\sO^{\circ}\to \sO^{\circ}/\omega$ induces a morphism $$\alpha_{\omega}(X):\sO^{\circ}(X)\to \sO^{\circ}/\omega(X).$$
Even when $X$ is affinoid, $\alpha_{\omega}(X)$ is not surjective in general.

\begin{exe}\label{exenotsurj}Let $\bk=\C((T)), \omega:=T$ and $X:=\sM(\bk\{x,y\}/(y^2-x^3-T^6)).$
Observe that $\sO^{\circ}(X)=\C[[T]][x,y]/(y^2-x^3-T^6)$ and $\sO^{\circ}(X)/\omega=\C[x,y]/(y^2-x^3).$  
Set $U_1:=\{z\in X|\,\, |x^3|\leq |y^2|\}$ and $U_2:=\{z\in X|\,\, |y^2|\leq |x^3|\}$. Both $U_1$ and $U_2$ are rational subdomain of $X.$
Define $F_1:=x^2/y\in O(U_1)$ and $F_2=y/x\in O(U_2)$.

We claim $\rho_{U_1}(F_1)\leq 1$ and $\rho_{U_2}(F_2)\leq 1.$
For every $z\in X$, if $R:=|x^3(z)|=|y^2(z)|$, then $0<R\leq 1, z\in U_1\cap U_2$ and $$|F_1(z)|=|(x^2/y)(z)|=|(y/x)(z)|=|F_2(z)|=R^{1/6}\leq 1.$$
For every $z\in X$, if $|x^3(z)|<|y^2(z)|$, then $z\in U_1\setminus U_2$, $|x(z)|<|T^2|$ and $|y(z)|=|T^3|$. So $|F_1(z)|=|(x^2/y)(z)|<|T|<1.$
For every $z\in X$, if $|x^3(z)|>|y^2(z)|$, then $z\in U_2\setminus U_1$, $|x(z)|=|T^2|$ and $|y(z)|<|T^3|$. So $|F_2(z)|=|(y/x)(z)|<|T|<1.$
Then we proved our claim.

For $z\in U_1\cap U_2$, we have $|x^3(z)|=|y^2(z)|\geq |T^6|.$ Then we have  $|xy(z)|\geq |T^5|.$
We have $F_1-F_2=x^2/y-y/x=(x^3-y^2)/xy=-T^6/xy.$
So we have $|(F_1-F_2)(z)|=|(T^6/xy)(z)|\leq |T^6|/|T^5|\leq |T|.$ It shows that $\alpha_{\omega}(F_1)|_{U_1\cap U_2}=\alpha_{\omega}(F_2)|_{U_1\cap U_2}.$
So we may glue $\alpha_{\omega}(F_1)$ and $\alpha_{\omega}(F_1)$ to a section $F\in \sO^{\circ}/\omega(X).$ 

We claim that $F$ is not contained in the image of $\alpha_{\omega}(X).$ Otherwise  $F=\alpha_{\omega}(X)(G)$ where $G\in \C[[T]][x,y]/(y^2-x^3-T^6).$
Define $U:=\{z\in X|\,\, |x(z)|=1\}.$ We have $U\subseteq U_1\cap U_2.$ Then on $U$, we get $$y/x= F|_U=G \mod \omega\in  \C[[T]][x,y]/(y^2-x^3-T^6)/(T)=\C[x,y]/(y^2-x^3)$$ which is a contradiction.

\end{exe}

The proof of Theorem \ref{thmisolim} relies on the following lemma.
\begin{lem}\label{lemdnsectionglo}Let $\omega$ be any element in $\bk$ satisfying $|\omega|\in (0,1).$ Then the morphism 
$$\alpha_{\omega}(\D^r):\sO^{\circ}(\D^r)\to \sO^{\circ}/\omega(\D^r)$$ is surjective for every $r\geq 0.$
\end{lem}

As pointed out by the referee of the first version of the paper, 
the rigid space version of this lemma is \cite[Theorem 1]{Bartenwerfer1978}.
The rigid space version means that we only consider strict covers.
Berkovich's trick of tensoring a suitable field $K_{R_1,\dots, R_s}$ \cite[Page 22]{Berkovich1990} allows us to reduce the Lemma \ref{lemdnsectionglo} to its rigid space version.
In this paper we give a different (and more direct) proof in the language of Berkovich spaces. 

\medskip

The following is a simple application of Lemma \ref{lemdnsectionglo}.
\begin{cor}\label{coraffineline}Let $\omega$ be any element in $\bk$ satisfying $|\omega|\in (0,1).$ Then we have $\sO^{\circ}/\omega(\A^{r,\an})=\bk^{\circ}/\omega$ for $r\geq 0.$
\end{cor}

\subsection*{Acknowledgement}
I would like to thank Bernard Le Stum, Xinyi Yuan, J\'er\^ome Poineau for useful discussions.
I thank Antoine Ducros for telling me that we need the condition ``generically quasi-smooth" in Lemma \ref{lemnoethergenerically}.
I thank Ruochuan Liu for sharing me his idea of an alternative proof of Theorem 1.1 in mixed characteristic using relative p-adic Hodge theory.
I thank the anonymous referee for the first version of the paper who point out that the rigid space version of our Theorem \ref{thmisolim} has a simplified  proof using \cite[Theorem 2]{Bartenwerfer1978}.

\section{A Noether normalization lemma}
Assume that the norm on $\bk$ is non-trivial.
Let $X=\sM(A)$ be a strict $\bk$-affinoid space.
We say that $X$ is \emph{generically quasi-smooth} if there is a Zariski dense open subset $U$ of $X$ such that $U$ is quasi-smooth.
\begin{rem}\label{remperfectgsmooth}If $\bk$ is perfect and $X$ is reduced, then $X$ is generically quasi-smooth.
\end{rem}

We need the following version of Noether normalization lemma for the proof of Theorem \ref{thmisolim}.

\begin{lem}\label{lemnoethergenerically}Assume that $X$ is irreducible and generically quasi-smooth of dimension $r$. Then there exists a finite morphism 
$$\pi: X\to \D^r=\sM(\bk\{T_1,\dots,T_r\}),$$ such that $\pi^*: \sO(\D^r)=\bk\{T_1,\dots,T_r\}\hookrightarrow \sO(X)$ is injective and
the field extension $\pi^*: \Frac \sO(\D^r)\hookrightarrow \Frac(\sO(X))$ is separable.
\end{lem}

\proof[Proof of Lemma \ref{lemnoethergenerically}]
By \cite[6.1.2. Theorem 1.(i)]{BGR}, there exists an admissible surjection $F: \bk\{U_1, \dots, U_n\}\twoheadrightarrow  A$ whose restriction $G:=F|_{\bk\{U_1, \dots, U_r\}}$ is finite and injective.
Since $X$ is generically quasi-smooth, there exists a rigid point $x\in X$ such that $X$ is quasi-smooth at $x$.
Then there exists an affinoid subdomain $W$ of $\D^n=\sM(\bk\{U_1, \dots, U_n\})$ containing $x$ such that $X\cap W$ is defined by the ideal $(g_1,\dots, g_{n-r})$ and moreover the Jacobian matrix at $x$
$J(x):=(\partial_{U_j}g_i(x))_{1\leq i\leq n-r, 1\leq j\leq n}\in M_{n\times (n-r)}(\sH(x))$ is of full rank. For $j=1,\dots, n$, set $J_j:=(\partial_{U_j}g_1(x),\dots, \partial_{U_j}g_{n-r}(x))$.
Then there are $1\leq j_1< \dots <j_{n-r}\leq n$ such that the vectors $J_{j_1},\dots, J_{j_{n-r}}$ are linearly independent.
Write $\{1,\dots,n\}\setminus \{j_1,\dots,j_{n-r}\}=\{s_1,\dots,s_r\}$ where $s_1<\dots< s_r.$
For $a\in \bk$, denote by $\Phi_a:\D^n\to \D^r$ sending $(U_1,\dots,U_n)$ to $(U_1+aU_{s_1},\dots, U_r+aU_{s_r})$ and set $\pi_a:=\Phi_a|_{X}.$
Then there exists $R>0$ such that when $0<|a|<R$, $(d\pi_a)(x)$ is invertible and the reduction $\widetilde{\pi_a^*}:\widetilde{\bk\{T_1,\dots, T_r\}}\to \widetilde{A}$ equal to the reduction $\widetilde{G}:\widetilde{\bk\{T_1,\dots, T_r\}}\to \widetilde{A}$. Define $\pi:=\pi_a$ for some $a\in \bk$ satisfying $0<|a|<R$.
Since $G$ is finite, $\widetilde{G}=\widetilde{\pi^*}$ is finite \cite[6.3.5. Theorem 1]{BGR}. Then $\pi$ is finite \cite[6.3.5. Theorem 1]{BGR}.
Since $(d\pi)(x)$ is invertible, the field extension $\pi^*: \Frac \sO(\D^r)\hookrightarrow \Frac(\sO(X))$ is separable, which concludes the proof.
\endproof

\section{The piecewise analytic functions}
Let $X$ be a reduced $\bk$-analytic space.  To avoid some set theoretical problem, we fix a cardinal number $\kappa\geq  \max\{\aleph_0, 2^{\# X}\}.$

\medskip

A \emph{piecewise analytic function} on $X$ is defined to be a collection $\{(U_i,f_i), i\in I\}$ indexed by an ordinal number $I$ of cardinal at most $\kappa$, where $U_i, i\in I$ is a G-covering of $X$ (by analytic domains) and $f_i\in \sO(U_i).$ 
Denote by $\sP(X)$ the set of piecewise analytic functions on $X$.
In particular, for $F=\{(U_i,f_i), i\in I\}, G=\{(V_i,g_i), i\in J\}$, we say that $F=G$ if $I=J$ and for every $i\in I=J$, $U_i=V_i$ and $f_i=g_i.$

For $F=\{(U_i,f_i), i\in I\}$, denote by $\sU_F:=\{U_i,i\in I\}$ the covering associated to $F.$
For $F=\{(U_i,f_i), i\in I\}, G=\{(V_i,g_i), j\in J\}$, we define $$F+G=\{(U_i\cap V_j, f_i+g_j), (i,j)\in I\times J\}$$ and 
$$F\times G=\{(U_i\cap V_j, f_ig_j), (i,j)\in I\times J\}.$$
Here $I\times J$ is the multiplication of ordinal numbers i.e. the order type of $I\times J$ is the lexicographical order.

\begin{rem}
When $\#I=1$ or $\#J=1$, we have $F\times G=G\times F.$
In general, $F\times G\neq G\times F$. Because usually the ordering on $I\times J$ and $J\times I$ are different.
\end{rem}

For $F,G,H\in \sP(X)$, we have 
$$F+G=G+F,  (F+G)+H=F+(G+H) \text{ and } (F\times G)\times H=F\times (G\times H).$$
In general $F\times (G+H)\neq F\times G+F\times H.$

\medskip

We have an embedding $\sO(X)\hookrightarrow \sP(X)$, sending $h$ to $\{(X,h)\}.$ This embedding preserves the addition and the  multiplication.
For every $F\in \sO(X), G,H\in \sP(X)$, we have $$F\times (G+H)= F\times G+F\times H.$$
We note that, for $F=\{(U_i,f_i), i\in I\}$, we have 
$$F-F=\{(U_i\cap U_j,f_i-f_j), (i,j)\in I\times J\}$$
and $$0\times F=\{(U_i,0), i\in I\}$$
are not necessarily equal to $0$.

\medskip

For $F=\{(U_i,f_i), i\in I\}\in \sP(X),$ define 
$$\rho(F)=\sup_{i\in I}\rho(f_i)\in [0,+\infty]$$
and $$\Delta(F):=\rho(F-F),$$ 
where $\rho$ is the spectral norm.
We may view $F$ as a multi-valued function on $X.$
For every $x\in X$, the value of $F$ at $x$ is denoted by $F(x):=\{f_i\in \sH(x)|\,\, x\in U_i\}.$
Then we have $$\rho(F)=\sup_{x\in X}\{|f_i(x)||\,\, x\in U_i\}$$ and $$\Delta(F)=\sup_{x\in X}\{|(f_i-f_j)(x)||\,\, x\in U_i\cap U_j\}.$$

For $F,G\in  \sP(X), a\in \bk$, we have 
$$\rho(F+G)\leq \max\{\rho(f),\rho(g)\},\rho(FG)\leq \rho(F)\rho(G);$$
$$\Delta(F+G)\leq \max\{\Delta(F), \Delta(G)\}, \Delta(FG)\leq \max\{\rho(F)\Delta(G),\rho(G)\Delta(F)\}.$$
and  $$\rho(aF)=|a|\rho(F), \Delta(aF)=|a|\Delta(F).$$

\medskip

For $s\in [0,\infty]$, denote by $\sP^s(X)$ (resp. $\sP_s(X)$ ) the set of $F\in \sP(X)$ satisfying $\rho(F)\leq s$ (resp. $\Delta(F)\leq s.$)
Set $\sP^s_t(X):=\sP^s(X)\cap \sP_t(X).$
They are stable under addition. It is easy to see the following properties
\begin{points}
\item $\sP(X)=\sP^{\infty}(X)$;
\item $\sP(X)^s_t$ is increasing for both $s$ and $t$;
\item $\sP^s_t(X)=\sP^{s}(X)$ if $s\leq t;$
\item $\sP^s(X)\subseteq \sP_s(X);$
\item for every $s_1,s_2,t_1,t_2\in [0,\infty]$, $\sP^{s_1}_{t_1}(X)\times \sP^{s_2}_{t_2}(X)\subseteq \sP^{s_1s_2}_{\max\{s_1t_2,s_2t_1\}}(X).$
\end{points}
For $F\in \sP(X)$,  $\Delta(F)=0$ if and only if  $F$ takes form $\{(U_i, f|_{U_i}), i\in I\}$ where $f\in \sO(X);$  if and only if there exists $f\in \sO(X)$ such that 
$\rho(F-f)=0.$

\medskip

Let $Y$ be a reduced $\bk$-analytic space satisfying $2^{\# Y}\leq \kappa$.
Let $f:Y\to X$ be a morphism. Then we have the pull back  map
$$f^*: \sP(X)\to \sP(Y)$$ sending $F:=\{(U_i,f_i), i\in I\}$ to $f^*F:=\{(f^{-1}(U_i),f^*(f_i)), i\in I\}.$
This map preserves the addition and the multiplication.
The restriction of $f^*$ on $\sO(X)$ is exactly the usual pull back of regular functions.
We have $\rho(f^*(F))\leq \rho(F)$ and $\Delta(f^*F)\leq \Delta(F).$

In particular, for any analytic subdomain $U$ of $X$, we may define the restriction map $\sP(X)\to \sP(U)$ by applying the pull back map for the inclusion.
Observe that, if $V_i, i\in I$ is a covering of analytic subdomains of $X$ and $F\in \sP(X)$, we have 
$$\rho(F)=\max\{\rho(F|_{V_i}), i=1,\dots, m\}.$$

\medskip

Let $I$ be an ordinal number with cardinal at most $\kappa$ and $V_i, i\in I$ be a G-covering of $X$ by analytic subdomains and $F_i=\{(V^i_j, f^i_j), j\in J_i\}\in \sP(V_i)$, we define 
$$\vee_{i\in I}F_i:=\{(V^i_j, f^i_j), (i,j)\in \sqcup_{s\in I}J_s\}\in \sP(X).$$
Here we write $(i,j)\in \sqcup_{s\in I}J_s$ for $i\in I$ and $j\in J_i$; and 
we think $\sqcup_{i\in I}J_i$ as a well ordinal set with the following ordering: $(a,b)\leq (c,d)$ if $a< c$ or $a=c,b\leq d$ and identify it with its order type.

We have 
$$\rho(\vee_{i\in I}F_i)=\sup_{i\in I} \rho(F_i)$$ and 
$$\Delta(\vee_{i\in I}F_i)\leq \sup\{\rho(F_i|_{V_i\cap V_j}-F_j|_{V_i\cap V_j}), i,j\in I\}.$$

\medskip

Let $\sU=\{U_i,i\in I\}$ be a G-covering of $X$ by analytic subdomains.  
A refinement of $\sU$ is a G-covering $\sV=\{V_j, j\in J\}$ with a map $\sigma:I\to J$ such that 
$V_i\subseteq U_{\sigma(i)}.$
For $F=\{(U_i,f_i), i\in I\}\in \sP(X)$, if $J$ is an ordinal number of cardinal at most $\kappa$, $\alpha:=(\sV=\{V_j, j\in J\}, \sigma)$ is a refinement of $\sU_F$. We define $F_{\alpha}\in \sP(X)$ by 
$$F_{\alpha}:=\{(V_i, f_{\sigma(i)}|_{V_i}), j\in J\}.$$
We have $$\rho(F_{\alpha})\leq \rho(F), \Delta(F_{\alpha})\leq \Delta(F) \text{ and } \rho(F_{\alpha}-F)\leq \Delta(F) .$$

\medskip
Let $F_n, n\geq 0$ be a sequence in $\sP(X).$  We say that $F_n, n\geq 0$ is a \emph{Cauchy sequence} if for every $\epsilon>0$, there exists $N\geq 0$ such that for every $n,m\geq N$, we have $\rho(F_n-F_m)\leq \epsilon.$ In fact, $F_n, n\geq 0$ in $\sP(X)$ is Cauchy if and only if  $\lim\limits_{n\to \infty}\Delta(F_n)=0$ and $\lim\limits_{n\to \infty}\rho(F_n-F_{n+1})=0.$ 
For $f\in \sO(X)\subseteq \sP(X)$, we say that $\lim\limits_{n\to \infty}F_n=f$ if $\rho(F_n-f)\to 0$ as $n\to \infty.$
It is clear that if $\lim\limits_{n\to \infty}F_n=f$, then we have 
\begin{equation}\label{rhofnf}
\rho(F_n-f)=\lim_{m\to \infty}\rho(F_n-F_m).
\end{equation}
Two Cauchy sequences $F_n, n\geq 0$ and $G_n, n\geq 0$ are defined to be equivalent if the 
$\lim\limits_{n\to \infty}(F_n-G_n)=0.$
Denote by $\hat{\sO}(X)$ the space of equivalent classes of Cauchy sequences in $\sP(X).$ 
For every Cauchy sequences $F_n, n\geq 0$ in $\sP(X)$, we denote by $[F_n]$ its image in $\hat{\sO}(X).$
Then we  have a natural inclusion $\sO(X)\subseteq \hat{\sO}(X)$ sending $f\in \sO(X)$ to the image of the constant sequence $F_n:=f, n\geq 0.$
We identify $\sO(X)$ with its image in $\hat{\sO}(X).$

It is easy to see that the addition and multiplication in $\sP(X)$ induce addition and multiplication on $\hat{\sO}(X)$ such that for every Cauchy sequences $F_n, G_n, n\geq 0$ in $\sP(X)$, we have 
$$[F_n]+[G_n]=[F_n+G_n] \text{ and } [F_n][G_n]=[F_nG_n].$$ The norm $\rho$ induces a norm on $\sP(X)$ (which is still denoted by $\rho$), satisfying 
$$\rho([F_n])=\lim_{n\to \infty}\rho(F_n).$$ 
Observe that $\rho([F_n])=0$ if and only if $[F_n]=0.$

We may view $[F_n]$ as a function on $X.$
Indeed, for every point $x\in X$, there exists a unique limit point $[F_n](x)$ for the sequences of subset $F_n(x)$ in $\sH(x)$. This point does not depend on the choice of $F_n, n\geq 0$ in its equivalent class. We say that $[F_n](x)$ is the value of $[F_n]$ at $x.$
Then we have 
\begin{equation}\label{eqdxnorm}
\rho([F_n])=\sup_{x\in X}|[F_n](x)|.
\end{equation}
Since $\rho([F_n]-[G_n])=\sup_{x\in X}|[F_n-G_n](x)|=\sup_{x\in X}|[F_n](x)-[G_n](x)|,$ $[F_n]=[G_n]$ if and only if they takes the same value at every point in $X.$

The restriction of $\rho$ on $\sO(X)$ is the spectral norm on $\sO(X)$. When $X$ is compact, $\rho$ takes finite values on $\hat{\sO}(X)$ and it makes $\hat{\sO}(X)$ to be a Banach $\bk$-algebra.

\medskip

The spaces $\sP(X), \sP^s(X),\sP_t(X),\sP^s_t(X)$ depends on the choice of the cardinal number $\kappa$. However, we have the following result.
\begin{pro}\label{prodxnotkap}The space $\hat{\sO}(X)$ does not depend on the choice of $\kappa$.
\end{pro}
\proof[Proof of Proposition \ref{prodxnotkap}]
Let $\beta$ be a cardinal number at least $\kappa.$ Denote by $\sP(X)'$ and $\hat{\sO}(X)'$ the spaces constructed in the same way as  $\sP(X)$ and $\hat{\sO}(X)$ except replacing $\kappa$ by $\beta.$ 
We have a natural inclusion $\sP(X)\subseteq \sP(X)'$. It induces a natural morphism $\hat{\sO}(X)\to \hat{\sO}(X)'$. Easy to see that this morphism is well defined. Moreover, it is a morphism of $\bk$-algebra and it preserves the spectral norm $\rho$. In particular, this morphism is injective. Then we may view $\hat{\sO}(X)$ as a subspace of $\hat{\sO}(X)'$. 
\begin{lem}\label{lemappkapp}For every $F'\in \sP(X)'$, there exists $F\in \sP(X)$ such that $\rho(F-F')\leq \Delta(F').$
\end{lem}
For every $[F_n']\in \hat{\sO}(X)'$, pick $F_n\in \sP(X), n\geq 0$ such that $\rho(F_n-F'_n)\leq \Delta(F'_n).$ 
Since $\Delta(F'_n)\to 0$ as $n\to \infty$, $F_n$ is a Cauchy sequence and $[F_n']=[F_n]\in \hat{\sO}(X).$ This concludes the proof.
\endproof
\proof[Proof of Lemma \ref{lemappkapp}]
Write $F'=\{(U_i,F_i),i\in I\}$.
Since the cardinal of the set of analytic subdomains of $X$ is at most $2^{\#X}\leq \kappa$, 
there exists a subset $J$ of $I$ of cardinal at most $\kappa$
such that for every $i\in I$, there exists $j\in J$ such that $U_i=U_j.$
We note that $J$ is a well-order set. We identify it with its order type. Define $F:=\{(U_j,F_j),i\in J\}\in \sP(X).$ 
We get $\rho(F-F')\leq \Delta(F'),$ which concludes the proof.
\endproof

For every $s>0$, denote by $\hat{\sO}^s(X)$ the set of $F\in \hat{\sO}(X)$ satisfying $\rho(F)\leq s.$
For every $F\in \hat{\sO}^s(X)$, we may write is as $F=[F_n]$ where $F_n$ is a Cauchy sequence in $\sP^s(X).$
Indeed, we may write $F=[F'_n]$ for some Cauchy sequence in $\sP(X).$
Then there exists $N\geq 0$, such that $\rho(F'_n)\leq s$ for $n\geq N.$ Then we have $F=[F'_{N+n}]$ and $F'_{N+n}\in \sP^s(X)$ for $n\geq 0.$
\begin{rem}
If we define $\widehat{\sO^s}(X)$ to be the equivalent classes of Cauchy sequences in $\sP^s(X)$, the above argument just shows that $\widehat{\sO^s}(X)=\hat{\sO}^s(X).$
\end{rem}

\medskip

\subsection{Sheaf properties}
\begin{pro}\label{procechsheaf}
The following presheaves 
$$\sP \text{ (resp. } \sP^s,\sP_t,\sP^s_t,\hat{\sO}, \hat{\sO}^s \text{)}: U\mapsto \sP(U) \text{ (resp. } \sP^s(U),\sP_t(U),\sP^s_t(U),\hat{\sO}(U), \hat{\sO}^s(U)\text{)}$$
are sheaves on $X$.
\end{pro}

\begin{rem}When $X$ is a complex manifold and $\sO$ is the sheaf of homeomorphisms, we may define the presheaf $\hat{\sO}$ in a same way. But usually it is not a sheaf and its sheafification is the sheaf of continuous functions.
\end{rem}

\proof[Proof of Proposition \ref{procechsheaf}]
We first prove it for $\sP$. We need to show that for any analytic subdomain $U$ of $X$,  every $G$-covering $V_j,j\in J$ of $U$, the first arrow in the following diagram is an equalizer:
$$\sP(U) \rightarrow \prod_{i\in J} \sP(V_i) {{{} \atop \longrightarrow}\atop{\longrightarrow \atop {}}} \prod_{i, j\in J} \sP(V_i \cap V_j).$$
We first show that the first arrow is injective. 
Let $F=\{(U_i, F_i), i\in I_1)\}, G=\{(W_i, G_i), i\in I_2)\}$ be two sections of $\sP(U),$ such that 
for every $j\in J$, $F|_{V_j}=G|_{V_j}.$ Then we get $I_1=I_2$. Moreover, for every $i\in I_1=I_2$, we have $U_i\cap V_j=W_i\cap V_j$ and $F_i|_{U_i\cap V_j}=G_i|_{W_i\cap V_j}.$
It follows that $U_i=W_i$ and $F_i=G_i$ for all $i\in I_1=I_2$. Then we get $F=G.$
Now assume that $(F_j, j\in J)\in \prod_{j} \sP(V_j)$ is contained in the equalizer. We only need to show that it is the image of some section of $\sP(U).$
Write $F_j=\{(U^j_i, F^j_i), i\in I_j\}.$ For $j,k\in J$, we have $F_j|_{V_j\cap V_k}=F_k|_{V_j\cap V_k}.$ It follows that $I:=I_j$ are the same  for all $j\in J.$ 
For every $i\in I$, define $U_i:=\cup_{j\in J}U^j_i.$ 
Since for every $i\in I, j,k\in J$, we have $U^j_i\cap V_j\cap V_k=U^k_i\cap V_j\cap V_k$, then for every $j\in J$, $U^j_i=U_i\cap V_j.$
Since for every $ j,k\in J$ $F^j_i|_{U_i\cap V_j\cap V_k}=F^k_i|_{U_i\cap V_j\cap V_k}$, there exists $F_i\in \sO(U_i)$ such that $F^j_i=F_i|_{U_i\cap V_j}.$
It follows that $(F_j, j\in J)$ is the image of $\{(F_i,U_i),i\in I\}\in \sP(U).$
We conclude the proof for $\sP.$
The proofs for $\sP^s,\sP_t,\sP^s_t$ are the same as the proof for $\sP$.

\medskip

We only need to prove it for $\hat{\sO}.$
We need to show that for any analytic subdomain $U$ of $X$,  every $G$-covering $V_j,j\in J$ of $U$, the first arrow in the following diagram is an equalizer:
$$\hat{\sO}(U) \rightarrow \prod_{i\in J} \hat{\sO}(V_i) {{{} \atop \longrightarrow}\atop{\longrightarrow \atop {}}} \prod_{i, j\in J} \hat{\sO}(V_i \cap V_j).$$
By Equality \ref{eqdxnorm}, the first arrow is injective. By Proposition \ref{prodxnotkap}, we may assume that $\kappa\geq \#J.$ We may endow $J$ a well-order and assume it to be an ordinal number.
Assume that $([F_n^j], j\in J)$ is an element in the equalizer in $\prod_{i\in J} \hat{\sO}(V_i).$
Define $F_n:=\vee_{j\in J}F_n^j\in \sP(X).$ It is easy to check that $F_n, n\geq 0$ is a Cauchy sequence and $[F_n]|_{V_j}=[F^j_n]|_{V_j}.$
We conclude the proof.
\endproof

\bigskip

Let $\omega$ be an element in $\bk$ satisfying $0<|\omega|<1.$ 
For every element $F=\{(U_i,F_i), i\in I\}\in \sP^1_{|\omega|}(X)$, we denote by $\widetilde{F}$ the section of $\sO_X^{\circ}/\omega(X)$ which is glued by $\widetilde{F_i}\in \sO_X^{\circ}/\omega(U_i)$ where 
$\widetilde{F_i}$ is the image of $F_i$ in $\sO_X^{\circ}(U_i)/\omega.$ 
Denote by $\psi_X: \sP^1_{|\omega|}(X)\to \sO_X^{\circ}/\omega(X)$ the map defined by $F\mapsto \widetilde{F}.$ 
Sometimes we write $\psi^{\omega}_X$ for $\psi_X$ to emphasize $\omega.$
We note that for $F_1,F_2\in \sP^1_{|\omega|}(X)$, $\psi_X(F_1)=\psi_X(F_2)$ if and only if $\rho(F_1-F_2)\leq |\omega|$.
We have the following result.
\begin{pro}\label{prooqomega}
The map $\psi_X:\sP^1_{|\omega|}(X)\to \sO^{\circ}/\omega(X)$ (resp. $\sP_{|\omega|}(X)\to \sO/\omega(X)$) is surjective.
\end{pro}
\proof[Proof of Proposition \ref{prooqomega}]
The $\sP^1_{|\omega|}$ part of this proposition implies the $\sP_{|\omega|}$ part. So we only need to proof the $\sP^1_{|\omega|}$ part.

Let $\widetilde{F}$ be a section in $\sO/\omega(X)$. Since the quotient morphism $\phi:\sO^{\circ}\to \sO^{\circ}/\omega(X)$ is surjective, there exists a $G$-covering 
$V_j, j\in J$ of $X$ and $F_j\in \sO^{\circ}(V_j), j\in J$ such that $\phi|_{V_j}(F_j)=\widetilde{F}|_{V_j}.$
Denote by $\beta$ a cardinal number at least $\max\{\#J,\kappa\}.$
Denote by $\sP(X)',\sP(X)'^1_{\omega}$ the spaces constructed in the same way as  $\sP(X), \sP(X)^1_{\omega}$ except replacing $\kappa$ by $\beta.$ 
View $\sP(X), \sP(X)^1_{\omega}$ as subspaces of $\sP(X)',\sP(X)'^1_{\omega}$. Define the morphism $\psi'_X:\sP'^1_{|\omega|}(X)\to \sO^{\circ}/\omega(X)$ in the same way as $\psi_X.$
We have $\psi'_X|_{\sP(X)^1_{\omega}}=\psi_X.$

After giving a well-order on $J$ and identifying it with its order type, we may assume that $J$ is an ordinal number of cardinal at most $\beta.$
Then $F':=\{(V_j,F_j),j\in J\}$ defines an element in ${\sP'^1_{|\omega|}}(X)$ satisfying 
$\psi_X'(F')=\widetilde{F}.$ By Lemma \ref{lemappkapp}, there exists $F\in \sP^1_{\omega}(X)$ such that $\rho(F-F')\leq \Delta(F')\leq |\omega|.$
Then we have $\psi_X(F)=\psi_X'(F)=\psi_X'(F')=\widetilde{F}$, which concludes the proof.
\endproof

For every analytic subdomain $U$ of $X$, $[F_n]\in \hat{\sO}^{1}(U)$, there exists $N\geq 0$ such that 
$F_n\in \sP_{|\omega|}^1(U)$ for $n\geq N$ and $\psi_X^{\omega}(F_n)=\psi_X^{\omega}(F_N).$
The value $\psi_X^{\omega}(F):=\psi_X^{\omega}(F_n)$ does not depend on the choice of $N$ and the sequence $F_n$.
Denote by $\beta_U^n: \hat{\sO}^1(U)\to \sO^{\circ}/\omega^n(U)$ the morphism 
sending $F$ to $\psi_X^{\omega^n}(F).$ They define a morphism of 
sheaves $\beta^n: \hat{\sO}^1\to \sO^{\circ}/\omega^n$.  Moreover $\beta_n$ induces a morphism of sheaves 
$\beta:=\lim\limits_{\leftarrow}\beta_n: \hat{\sO}^1\to \lim\limits_{\leftarrow}\sO^{\circ}/\omega^n.$

\begin{pro}\label{prodxlimo} The morphism $\beta: \hat{\sO}^1\to \lim\limits_{\leftarrow}\sO^{\circ}/\omega^n$ (resp. $\hat{\sO}\to \lim\limits_{\leftarrow}\sO/\omega^n\sO^{\circ}$) is an isomorphism.
\end{pro}
\proof[Proof of Proposition \ref{prodxlimo}]
The $\hat{\sO}^1$ part of this proposition implies the $\hat{\sO}$ part. So we only need to proof the $\hat{\sO}^1$ part.

We only need to show that for affinoid subdomain $U$, $\beta_U:\hat{\sO}^1(U)\to \lim\limits_{\leftarrow}\sO^{\circ}/\omega^n(U)$ is an isomorphism.

For $F=[F_n]\in \hat{\sO}^1(U)$, we may assume that $\rho(F_n)\leq 1.$ If $F\neq 0$, there exists $l\geq 0$ such that $\rho(F)>|\omega^l|.$
There exists $N\geq 0$, such that for $n\geq N$, $\psi_X^{\omega^l}(F_n)=\psi_X^{\omega^l}(F_N).$  
Denote by $q_l: \lim\limits_{\leftarrow}\sO^{\circ}/\omega^n(U)\to \sO^{\circ}/\omega^l(U)$ the quotient morphism, we have 
$$q_l(\beta(F))=\beta_n(F)=\psi_X^{\omega^l}(F_N)\neq 0.$$ This shows that $\beta$ is injective.

For every $(g_n)_{n\geq 0}\in \lim\limits_{\leftarrow}(\sO^{\circ}/\omega^n(U))=\lim\limits_{\leftarrow}\sO^{\circ}/\omega^n(U)$, 
by Proposition \ref{prooqomega}, we may write $g_n=\psi_U^{\omega^n}(G_n)$ for some $G_n\in \sP_{|\omega^n|}^1(U).$
We may check that $G_n, n\geq 0$ is a Cauchy sequence in $\sP^1(U)$ and $\beta_U([G_n])=(g_n)_{n\geq 0}.$
Then $\beta$ is surjective, which concludes the proof.
\endproof

\subsection{An induced continuous $\R$-valued function}\label{subsectionphipb}
For $b>0$, $F=\{(U_i, F_i), i\in I\}\in \sP_b(X)$, we define a function $\phi_F^b: X\to \R$ on $X$ as follows:
for every $x\in X$, $$\phi_F^b(x):=\max\{\log(|F_i(x)|),\log (b)\},$$ where $x\in U_i.$
It is well defined, because for $i,j\in I$ satisfying $x\in U_i\cap U_j$, we have $|(F_i-F_j)(x)|\leq b$, then $\max\{\log(|F_i(x)|),\log |b|\}=\max\{\log(|F_j(x)|),\log |b|\}$.
We note that for $F'\in \sP_{b}(X)$, if $ \rho(F-F')\leq b$, then we have $\phi_F^b=\phi_{F'}^b$.

\begin{pro}\label{procontinuityphif}The function $\phi_F^b$ is continuous on $X$.
\end{pro}

\proof[Proof of Proposition \ref{procontinuityphif}]
For every $i\in I$, $\phi_F^b|_{U_i}=\max\{\log(|F_i|_{U_i}|),\log (b)\}$ is continuous. We conclude the proof by the following lemma.
\begin{lem}\label{lemcontixonlyneedaff}
Let $U_i, i\in I$ be an affinoid $G$-covering of $X$. 
Then a function $\phi: X\to \R$ is continuous if and only if for every $i\in I$, $\phi|_{U_i}: U_i\to \R$ is continuous.
\end{lem}
\endproof
\proof[Proof of Lemma \ref{lemcontixonlyneedaff}]
It is clear that if $\phi$ is continuous, $\phi|_{U_i}$ is continuous for every $i\in I.$

Now assume that $\phi|_{U_i}$ is continuous for every $i\in I.$
Let $K$ be a closed subset of $\R$, we only need to show that $\phi^{-1}(K)$ is closed in $X.$
Let $x$ be a point in $X\setminus \phi^{-1}(K)$, since $U_i, i\in I$ is a $G$-covering of $X$, there is a finite subset $I_x$ of $I$ such that $$x\in (\cap _{i\in I_x}U_i)\cup ((\cup_{i\in I_x}U_i)^{\circ}).$$
Then $(\cup_{i\in I_x}U_i)\cap  \phi^{-1}(K)=\cup_{i\in I_x}(\phi|_{U_i})^{-1}(K)$ is closed in $X$. So $\phi^{-1}(K)$ contains the open neighborhood $((\cup_{i\in I_x}U_i)^{\circ})\setminus \phi^{-1}(K)$ of $x$, which concludes the proof.
\endproof

\section{Proof of Lemma \ref{lemdnsectionglo} and Corollary \ref{coraffineline}}
\proof[Proof of Lemma \ref{lemdnsectionglo}]
Denote by $X:=\D^r=\sM(\bk\{T_1,\dots,T_r\}).$ Set $b:=|\omega|.$
When $r=0$, there is nothing to prove. Now we assume that $r\geq 1.$

By Proposition \ref{prooqomega}, we only need to show that for every element $F=\{(U_i,F_i), i\in I\}\in \sP^1_b(X)$, there exists a section $G\in \sO^{\circ}(X)$ such that 
for every $i\in I$, $\rho(G|_{U_i}-F_i)\leq b.$

\subsection{Reduce to the finite covering by rational subdomains}\label{subsectionreducefinite}
Since for every $i\in I$, $U_i$ has a $G$-covering of (at most $\kappa$) rational subdomains, after a refinement, we may assume that $U_i, i\in I$ are rational subdomains.
In particular, all $U_i, i\in I$ are compact.

For every $x\in X$, there exists a finite subset $I_x\subseteq I$ such that $$x\in (\cup_{i\in I_x} U_i)^{\circ}\cap (\cap_{i\in I_x}U_i).$$
Since $X$ is compact, there exists a finite set $S\subseteq X$, such that $X\subseteq \cup_{x\in S} (\cup_{i\in I_x} U_i)^{\circ}.$
\begin{lem}\label{lemcriquasinet}The collection $U_i, i\in \cup_{x\in S}I_x$ is a $G$-covering of $X.$ 
\end{lem}
\proof[Proof of Lemma \ref{lemcriquasinet}]
For every $x\in X$, there exists $y\in S$ such that 
$x\in (\cup_{i\in I_y} U_i)^{\circ}.$ Let $I_1$ be the set of  $i\in I_y$ such that $x\in U_i.$
Then we get $x\in \cap_{i\in I_1}U_i$
and 
$$x\in (\cup_{i\in I_y} U_i)^{\circ}\setminus (\cup_{j\in I_y\setminus I_1}U_j)\subseteq \cup_{i\in I_1}U_i.$$
Since $(\cup_{i\in I_y} U_i)^{\circ}\setminus (\cup_{j\in I_y\setminus I_1}U_j)$ is open, we get 
$x\in (\cup_{i\in I_1}U_i)^{\circ}.$ This concludes the proof.
\endproof
After replacing $I$ by $\cup_{x\in S}I_x$, we may assume that $I$ is finite.

%

\subsection{The case $r=1$}\label{subsectionrone}
In this case, $X=\sM(\bk\{T\})$ is a rooted $\R$-tree. Denote by $\leq$ the natural partial ordering on $X$ by declaring that $x\leq y$ if and only if $|f(x)|\leq |f(y)|$ for all $f\in \bk\{T\}$.
The Gauss point $\xi: f\in \bk\{T\} \mapsto \rho(f)$ is the unique maximal element. 

For every element $H=\{(U_i',H_i),i\in I'\}\in \sP_{b}(X)$ and $c\geq b$, define $\phi_H^c: X\to \R$ as in
Section \ref{subsectionphipb}.

%

\begin{lem}\label{lemquasisubhar}For every $c\geq b$, the function $\phi_H^c$ is $\xi$-subharmonic on $X.$
\end{lem}
\proof[Proof of Lemma \ref{lemquasisubhar}]
We first treat the case where $c>b.$

By the argument in Section \ref{subsectionreducefinite}, we may assume that $I'$ is finite and $U'_i,i\in I$ are rational subdomains. 
Since $U_i',i\in I'$ are rational, $\bk(T)\cap \sO(U_i')$ is dense in $\sO(U_i').$
So we may assume that $H_i\in \bk(T)\cap \sO(U_i').$
Set $I_i:=\{j\in I'|\,\, U_i'\cap U_j'\neq \emptyset\}.$ 
Since $U_i', i\in I'$ are compact, there are open sets $V_i, i\in I'$ such that 
\begin{points}
\item $U_i'\subseteq V_i$;
\item $V_i\cap V_j=\emptyset$ if and only if $U_i'\cap U_j'=\emptyset.$
\end{points}
For $i\in I'$, define an open set $W_i:=\cap_{j\in I_i}\{x\in V_i|\,\,|H_i-H_j|<c\}.$ 
Since $U_i'\subseteq W_i$, $W_i, i\in I$ forms an open covering of $X.$
For every $x\in W_i$, we have $\phi_H^c=\max\{\log(|H_i|), \log(c)\},$ which is $\xi$-subharmonic.
Then $\phi_H^c$ is $\xi$-subharmonic.

Now we only need to treat the case $b=c$.
Since $\phi_H^b$ is the uniform limit of the sequence of $\xi$-subharmonic functions $\phi_H^{b+(1-b)/n}, n\geq 0$,  $\phi_H^b$ is $\xi$-subharmonic.
This concludes the proof.
\endproof

Since $U_i,i\in I$ are rational, $\bk(T)\cap \sO(U_i)$ is dense in $\sO(U_i).$
So we may assume that $F_i\in \bk(T)\cap \sO(U_i).$ We may assume that $\xi\in U_0$.

\begin{lem}\label{lemapprorootleone}For every $\epsilon>0$, there exists a rational function $P/Q\in \bk(T)\cap \sO(U_0), P,Q\in \bk[T]$ such that 
$\rho_0(F_0-P/Q)\leq \epsilon$ and the norms of all roots of $Q=0$ are at most $1$.
\end{lem}

\proof[Proof of Lemma \ref{lemapprorootleone}]
Write $F_0=P/Q$ where $P,Q\in \bk[T]$ and $Q\neq 0.$
We may assume that $Q$ is monic. We note that for every monic irreducible factor $Q'$ of $Q$, the norm of all roots of $Q'=0$ are the same.
So we may write $Q=Q_1Q_2$ where $Q_i,i=1,2$ are monic polynomials such that the norms of all roots of $Q_1=0$ are at most $1$ and all roots of $Q_2=0$ are strictly larger than $1$.
Since $1/Q_2\in \sO(\D)=\bk\{T\}$, after replacing $Q_2^{-1}$ by some polynomial in $\bk[T]$ which is close enough to $1/Q_2\in \bk\{T\}$, we may assume that $Q=Q_1$. In other words, we may assume that the norms of all roots of $Q=0$ are at most $1$. This concludes the proof.
\endproof

By Euclidean division, we may write $P=GQ+R$ where $G, R\in \bk[T]$ and $\deg R<\deg Q.$
Then we have $F_0=G+R/Q.$ Set $L:=F-G\in \sP_b(X).$
By Lemma \ref{lemquasisubhar}, $\phi_L^b$ is a $\xi$-subharmonic on $X.$
In particular, we have $$\sum_{v\in  {\Tan_{\xi}}X}D_v\phi_L^c\leq 0.$$
\begin{lem}\label{labelsmallb}We have $|R/Q(\xi)|\leq b.$
\end{lem}
\proof[Proof of Lemma \ref{lemquasisubhar}]
Assume that $|R/Q(\xi)|>b$.

For every $i\in I$ such that $\xi\in U_i$, $L_i=F_i-G$ takes form
$L_i=R/Q+\omega S$ where $S\in \bk(T)$ and $|S(\xi)|\leq 1.$
Then there exists an open neighborhood $V$ of $\xi$, such that for every $i\in I, x\in V\cap U_i$, we have
$|L_i(x)|=|R/Q(x)|>b.$
Then, for $v\in {\Tan_{\xi}}X$, we have 
$$D_v\phi_L^b=D_v(\log |R/Q|)=D_v(\log |R|)-D_v(\log |Q|).$$
We have 
$$\sum_{v\in  {\Tan_{\xi}}X}D_v(\log |R|)\geq -\deg R.$$
Since the norms of all roots of $Q=0$ are at most $1$, we have 
$$\sum_{v\in  {\Tan_{\xi}}X}D_v(\log |Q|)=-\deg Q.$$

It follows that 
$$\sum_{v\in  {\Tan_{\xi}}X}D_v\phi_L^b=\sum_{v\in  {\Tan_{\xi}}X}D_v(\log |R/Q|)\geq \deg Q-\deg R>0,$$
which is a contradiction.
\endproof
By Lemma \ref{lemquasisubhar},  we have 
$$ \phi_L^b(\xi)=\max\{\log |L_0(\xi)|,  \log (b)\}=\max\{\log |R/Q(\xi)|, \log (b)\}= \log (b).$$
Since $\phi_L^b$ is decreasing on $X$, we get $\phi_L^b=\log b$ on $X.$
It follows that $$\rho(F-G)=\rho(L)\leq b<1.$$
It follows that $G\in \sO^{\circ}_X(X)$ and $\rho(G|_{U_i}-F_i)\leq b$ for $i\in I.$ This concludes the proof when $r=1.$
\subsection{The general case}
Now we prove Lemma \ref{lemdnsectionglo} by induction on $r\geq 1.$
We may assume that $r\geq 2.$

Since $U_i,i\in I$ is rational, $\bk(T_1,T_2,\dots, T_r)\cap \sO(U_i)$ is dense in $\sO(U_i).$ 
We may assume that $F_i\in \bk(T_1,T_2,\dots, T_r)$ for all $i\in I.$

Consider the morphism $$\pi: X=\sM(\bk\{T_1,\dots,T_r\})\to \D:=\sM(\bk\{T_1\})$$ defined by $$\pi^*: T_1\mapsto T_1.$$
For every $x\in \D$, denote by $$X_x:=\sM(\sH (x)\{T_2,\dots,T_r\})\simeq \D_{\sH (x)}^{r-1}$$ the fiber of $\pi$ above $x.$
We have $F_x:=F|_{X_x}=\{(U_i\cap X_x,F_i|_{U_i\cap X_x}), i\in I\}\in \sP^1_b(X_x).$ 
For  every $c\geq b$, define a function $\psi_F^c: \D\to \R$ by 
$$x\mapsto \max\{\rho(F|_{X_x}),\log(c)\}.$$
Let $\xi$ be the Gauss point in $\D.$
\begin{lem}\label{lemsubharhigdim}For every $c\geq b$, the function $\psi_F^c$ is $\xi$-subharmonic on $\D.$
\end{lem}
\proof[Proof of Lemma \ref{lemsubharhigdim}]
We first treat the case where $c>b.$

By the induction hypotheses, for every $x\in \D$,  there exists $H_x\in \sO^{\circ}_{X_x}(X_x)=\sH (x)\{T_2,\dots,T_r\}$ such that 
$\rho(H_x-F_x)\leq b.$ There exists $g_x\in \bk[T_1]\setminus \{0\}$, $G_x\in \bk[T_1]_{(g_x)}[T_1,\dots, T_r]$
such that $g_x(x)\neq 0$ and  $\rho(H_x-G_x|_{X_x})\leq b.$
There exists a rational neighbourhood $W_x$ of $x$ in $V_x:=\{y\in \D|\,\, g_x(y)\neq 0\}$, such that 
\begin{points}
\item 
$\pi^{-1}(W_x)\cap U_i=\emptyset$ if and only if $X_x\cap U_i=\emptyset;$ 
\item for every $i\in I$, $\rho(G_x|_{\pi^{-1}(W_x)\cap U_i}-F_i|_{\pi^{-1}(W_x)\cap U_i})<c.$
\end{points}
Then for every $z\in \pi^{-1}(W_x)$, we have $$\max\{\log|G_x(z)|, \log (c)\}=\psi_F^c(z).$$
Since $\D$ is compact, there exists a finite set $\{x_0,\dots x_m\}\subseteq \D$ such that $$\cup_{j=0}^mW_{x_j}=\D.$$
Set $H:=\{(\pi^{-1}(W_{x_i}),H_i),i\in \{0,\dots,m\}\}\in \sP^1_c(X)$, where $H_i:=G_{x_i},$ we have $\psi^c(H)=\psi^c(F).$
Set $W_i:=W_{x_i}, i=0,\dots,m$. 
There exists $l\geq 1$, such that, for every $i=0,\dots,m$, 
$$H_i=\sum_{|J|\leq l}A_i^J(T_1)T^J$$
where $J=(j_2,\dots,j_r)\in \Z_{\geq 0}^{r-1}$ is a multi-index, $|J|=j_2+\dots+j_r$, $T^J:=T_2^{j_2}\dots T_r^{j_r}$ and $A_i^J(T_1)\in \bk(T_1)\cap \sO_{\D}(W_i).$
For $|J|\leq l,$ set $$H^J:=\{(W_i,A^J_i),i\in \{0,\dots,m\}\}\in \sP_c(\D).$$
Define $\phi_{H^J}^c$ as in Section \ref{subsectionrone}.
By Lemma \ref{lemquasisubhar},  $\phi_{H^J}^c$ is $\xi$-subharmonic.

For every $x\in W_i$, we have 
$$\max\{\log(\rho(H|_{X_x})),\log (c)\}=\max\{\log(\rho(H_i|_{X_x})),\log((c)\}$$
$$=\max_{|J|\leq l}\max\{\log |A_i^J(x)|,\log((c)\}=\max_{|J|\leq l}\phi_{H^J}^c.$$
It follows that $$\psi^c(F)=\psi^c(H)=\max_{|J|\leq l}\phi_{H^J}^c$$ is $\xi$-subharmonic.

Now we treat the case $c=b.$
Since $\psi_F^b$ is the uniform limit of the sequence of $\xi$-subharmonic functions $\phi_F^{b+(1-b)/n}, n\geq 0$,  $\phi_F^b$ is $\xi$-subharmonic.
This concludes the proof.
\endproof

By Lemma \ref{lemsubharhigdim}, we have $$\psi^b_F(\xi)=\sup_{x\in \D}\psi^b_F(x).$$
Let $\eta$ be the Gauss point of $X=\D^r$. It can be also viewed as the Gauss point of $X_{\xi}\simeq \D^{r-1}_{\sH(\eta)}.$
By the induction hypotheses, there exists $H\in \sO^{\circ}_{X_{\xi}}(X_{\xi})=\sH (\xi)\{T_2,\dots,T_r\}$ such that 
$\rho(H-F_{\xi})\leq b.$ 
It shows that 
\begin{equation}\label{eqmax}
\begin{split}
\sup_{x\in \D}\psi^b_F(x)=\psi^b_F(\xi)
=\max\{\log \rho(F|_{\xi}), \log (b)\}\\=\max\{\log \rho(H), \log (b)\}
=\max\{\log |H(\eta)|, \log(b)\}\\=\max\{\log |F(\eta)|, \log(b)\},
\end{split}
\end{equation}
where $|F(\eta)|:=\sup\{|F_i(\eta)||\,\, \eta\in U_i\}.$ This equation can be viewed as the maximal principle for elements in $\sP_b(\D^r).$

Since $\bk(T_1)$ is dense in $\sH (\xi)$,  we may assume that $H\in \bk(T_1)[T_2,\dots,T_r].$
Write $$H=G(T_1,\dots, T_r)+\sum_{|J|\leq l}A_J(T_1)/B(T_1)T^J$$ where $J=(j_2,\dots,j_r)\in \Z_{\geq 0}^{r-1}$ is a multi-index, $|J|=j_2+\dots+J_m$, $T^J:=T_2^{j_2}\dots T_r^{j_r}$, $G(T_1,\dots,T_r)\in \bk[T_1,\dots, T_n]$, $A_J(T_1), B(T_1)\in \bk[T_1]$. By Lemma \ref{lemapprorootleone}, we may assume that  the norms of all roots of $B(T_1)=0$ are at most one. By euclidean division, we may assume that $\deg A_J<\deg B.$
Set $L:=F-G\in \sP_b(X).$

We claim that $\max_{|J|\leq l}\{\rho(A_J(T_1)/B(T_1))\}\leq b.$ Otherwise there exists $J$, with $|K|\leq l$ and $\rho(A_K(T_1)/B(T_1))>c> b.$
There exists a rational neighbourhood $W$ of $\xi$, such that 
\begin{points}
\item 
$\pi^{-1}(W)\cap U_i=\emptyset$ if and only if $X_{\xi}\cap U_i=\emptyset;$ 
\item for every $i\in I$, $\rho(H|_{\pi^{-1}(W)\cap U_i}-F_i|_{\pi^{-1}(W)\cap U_i})<c;$
\item for every $x\in W$, $|A_K(x)/B(x)|>c.$
\end{points}
It follows that $R:=L\vee \{(\pi^{-1}(W), H-G)\}\in \sP_c(X).$
Apply Lemma \ref{lemsubharhigdim} for $R$, the function 
$\psi^c_R$ on $\D$ is $\xi$-subharmonic.
For $x\in W$, we have 
$$\psi^c_R(x)=\max\{\max_{|J|\leq l}\log |A_J(x)/B(x)|, \log (c)\}\geq \log |A_K(x)/B(x)|.$$
Then we have 
$$0\geq \sum_{v\in  {\Tan_{\xi}}\D}D_v\psi_R^c\geq \sum_{v\in  {\Tan_{\xi}}X}D_v(\log |A_K/B|)\geq\deg A_K-\deg B>0,$$
which is a contradiction. This concludes the claim.

Then we have 
$$\rho((F-G)|_{X_{\xi}})\leq \max\{\rho((H-G)|_{X_{\xi}}),\rho(F|_{X_{\xi}}-H)\}\leq b.$$
We note that $F-G\in \sP_b(X)$. Apply Equality \ref{eqmax} for it, we get $$\rho(F-G)\leq \rho((F-G)|_{X_{\xi}})\leq b.$$
This concludes the proof of Lemma \ref{lemdnsectionglo}.
\endproof
\subsection{Proof of Corollary \ref{coraffineline}}
When $r=0$, there is nothing to prove. Now assume that $r\geq 1.$
Set $b:=|\omega|.$
Let $F$ be a section in $\sO^{\circ}/\omega(\A^{r,\an}).$
For every $R\in (0,+\infty)\cap |\bk|$, denote by $\D^r(R):=\sM(\bk\{T_1/R,\dots,T_r/R\})$ the ball containing $0$ and of radius $R.$

Pick $\beta\in \bk$ with $c:=|\beta|>1.$

By Lemma \ref{lemdnsectionglo}, for every $n\geq 1$, there exists $F_n\in \bk^{\circ}\{T_1/c^n,\dots,T_r/c^n\}$ such that 
$$\alpha_n(F_n)=F|_{\D^r(c^n)},$$ where $\alpha_n:=\alpha_{\omega}(\D^r(c^n)).$
We may assume that $F_n\in \bk^{\circ}[T_1,\dots,T_r]$ and write $$F_n=\sum_{I}a^n_IT^I.$$
For $m\geq n$, we have $\rho(F_m|_{\D^r(c^n)}-F_n)\leq b.$
In other words, we have, for $m\geq n,$
\begin{equation}\label{eqinmain}
|a_I^m-a_I^n|\leq b/c^{n|I|}
\end{equation}
Since $\rho(F_n)\leq 1$ for $n\geq 1$, we have 
\begin{equation}\label{eqinmainto}
|a_I^n|\leq 1/c^{n|I|}.
\end{equation}

In particular, we have $|a_0^n|\leq 1$ and  $|a_0^n-a^0_0|\leq b$ for $n\geq 0.$
Denote by $\bar{a}$ the image of it in $\bk^{\circ}/\omega.$
After replacing $F$ by $F-\bar{a}$, we may assume that $a^n_0=0, n\geq 0.$
We only need to show that $F=0.$
We note that $F$ is represented by $H:=\{(\D^r(c^n), F_n), n\geq 1\}\in \sP^0_b(\A^{r,\an}).$
We only need to show that 
$\rho(H)\leq b.$

By Inequality \ref{eqinmainto},
for every $n\geq 1, m\geq n$, we have $$\max\{\rho(H|_{\D^r(c^n)}),b\}=\max\{\rho(F^m|_{\D^r(c^n)}),b\}$$
$$=\max\{\max_{I}|a_I^m||c|^n,b\}\leq \max\{\max_{|I|\geq 1}|c|^{n-m|I|},b\}=\max\{|c|^{n-m},b\}.$$
For $m$ large enough, we have $|c|^{n-m}<b$. Then we get $\rho(H|_{\D^r(c^n)})\leq b$ for $n\geq 0.$ It follows that $\rho(H)\leq b$, which conclude the proof.

\section{Proof of Theorem \ref{thmisolim}}\label{sectionproofthmmain}
Let $X$ be a reduced normal $\bk$-analytic space.
We need to show that for every affinoid subdomain $U$ of $X$, the morphism $\alpha_U: \sO^{\circ}(U)\to \lim\limits_{\leftarrow}\sO^{\circ}/\omega^n(U)$ is an isomorphism.

So we may assume that $X$ is an affinoid space and we only need to show that the morphism $\alpha_X: \sO^{\circ}(X)\to \lim\limits_{\leftarrow}\sO^{\circ}/\omega^n(X)$ is an isomorphism.

It is clear that $\alpha_X$ is injective.  So we only need to show that $\alpha_X$ is surjective.
By Proposition \ref{prodxlimo}, we identify $\lim\limits_{\leftarrow}\sO^{\circ}/\omega^n(X)$ with $\hat{\sO}^1(X)$ and identify $\sO^{\circ}(X)$ as its image in $\hat{\sO}^1(X).$
Let $F=[F_n]$ be an element in $\hat{\sO}^1(X)$. We may assume that $F_n\in \sP^1(X).$
Write $F_n=\{(U_i^n, F^n_i), i\in I^n\}.$
By the argument in Section \ref{subsectionreducefinite}, we may assume that $I^n$ are finite and $U_i^n$ are rational subdomains in $X$.
We may assume that for every $m\geq n\geq 0,$
$\rho(F_m-F_n)\leq 1/2^n.$

Since $U_i^n, n\geq 0, i\in I^n$ are rational domains, 
$\Frac(\sO(X))\cap \sO(U_i^n)$ are dense in $\sO(U_i^n)$.  We may assume that $F^n_i$ are contained in $\Frac(\sO(X))\cap \sO(U_i^n)$ for all $n\geq 0,i\in I^n.$

Let $X_1,\dots, X_m$ be the irreducible components of $X$. 
Since $X$ is normal, we have $X=\sqcup_{i=1}^m X_i.$
It is clear that $\sO^{\circ}(X)=\oplus_{i=1}^m\sO^{\circ}(X_i)$ and $\hat{\sO}^1(X)=\oplus_{i=1}^m\hat{\sO}^1(X_i).$
So we may assume that $X$ is irreducible.


\subsection{Reduce to the strict case}\label{subsectionredtostrict}
\begin{lem}\label{lemextrithm}
Assume that $R\in \R_{>0}\setminus \sqrt{|\bk^*|}.$  Let $K$ be the field $$K:=K_R:=\bk\{R^{-1}S_1,RS_2\}/(S_1S_2-1)$$
as in \cite[Page 21]{Berkovich1990}.
If Theorem \ref{thmisolim} holds for the $K$ analytic space $X_K:=X\widehat{\otimes}_{\bk}K$, then it holds for $X.$
\end{lem}
\proof[Proof of Lemma \ref{lemextrithm}]
By our assumption, there exists $G_K\in \sO^{\circ}(X)\otimes_{\bk}K$ such that for every $n\geq 0, i\in I^n$, $\rho(G_K|_{U^n_i\widehat{\otimes}_{\bk}K }-F^n_i)\leq 1/2^n.$

Write $G_K=\sum_{j\in \Z}G_jS_1^j,$ where $G_j\in \sO^{\circ}(X).$
For every $n\geq 0, i\in I^n$, we have $$\rho(\sum_{j\in \Z^*}G_j|_{U^n_i} S_1^j+(G_0-F_i^n))\leq 1/2^n.$$
By \cite[Section 2.1]{Berkovich1990}, we have $$\rho(\sum_{j\in \Z^*}G_j|_{U^n_i} S_1^j+(G_0-F_i^n))=\max\{\rho(G_0|_{U^n_i}-F^n_i), \max_{j\neq 0}\rho(G_j|_{U^n_i})R^j\}.$$
It follows that $\rho(G_0|_{U^n_i}-F^n_i)\leq 1/2^n, i\in I_n$, which concludes the proof.
\endproof
By this lemma, we may assume that $X$ is strict.

\subsection{Multiply by an analytic function}
It is clear that for $Q\in \sO^{\circ}(X)$ and $n\geq m\geq 0$, we have $$\rho(QF_n-QF_m)\leq 1/2^m.$$
So $QF_n, n\geq 0$ is a Cauchy sequence.  

\begin{lem}\label{lemmultiqanna}For $Q\in \sO^{\circ}(X)\setminus \{0\}$, if $QF\in \sO^{\circ}(X)\subseteq \hat{\sO}^1(X)$, then $F\in \sO^{\circ}(X).$
\end{lem}

\proof[Proof of Lemma \ref{lemmultiqanna}]
Set $U:=X\setminus \{x\in X|\,\, |Q(x)|=0\}.$
Since $QF\in \sO^{\circ}(X)$, $G:=(QF)/Q\in \sO(U).$
By Equality \ref{rhofnf}, 
\begin{equation}\label{equonufng}\rho(F_n|_U-G)\leq 1/2^n.
\end{equation}

For every $i\in I^n, x\in U_i$, we have $|F^n_i(x)|\leq 1.$
It follows that $|QF^n_i(x)|\leq |Q(x)|.$
Then for every $x\in U$, we have $|G(x)|\leq 1.$
Since $G\in \Frac(\sO(X))$ and $X$ is normal, by Riemann extension theorem \cite[3.3.14]{Berkovich1990}, we get $G\in \sO(X).$
Then $F_n-G\in \sP^1_{1/2^n}(X).$
Recall that for every $x\in U_i$, $$\phi^{1/2^n}_{F_n-G}(x)=\max\{\log |(F^n_i-G)(x)|, \log (1/2^n)\}.$$
By Inequality \ref{equonufng}, for $x\in U$, $\phi^{1/2^n}_{F_n-G}(x)=\log (1/2^n).$
By Proposition \ref{procontinuityphif},  $\phi^{1/2^n}_{F_n-G}: X\to \R$ is continuous. 
Since $X$ is irreducible, $U$ is dense in $X$. It follows that $\phi^{1/2^n}_{F_n-G}=\log (1/2^n).$ Then we have  
$\rho(F_n-G)\leq 1/2^n.$
This implies that $G=[F_n]$, which concludes the proof.
\endproof

\subsection{Proof of Theorem \ref{thmisolim}}
By Lemma \ref{lemnoethergenerically}, there exists a finite morphism 
$$\pi: X\to \D^r$$ where $r=\dim X$, such that $\pi^*: \sO(\D^r)=\bk\{T_1,\dots,T_r\}\hookrightarrow \sO(X)$ is injective and 
the field extension $\pi^*: \Frac \sO(\D^r)\hookrightarrow \Frac(\sO(X))$ to be separable.
View $\Frac \sO(\D^r)$ as a finitely dimensional vector space over $\Frac(\sO(X)).$
Denote by $\Tr: \Frac(\sO(X))\to \Frac \sO(\D^r)$ the trace morphism. The bilinear form 
$$\Tr(\cdot\times \cdot): \Frac(\sO(X))\times \Frac(\sO(X))\to \sO(\D^r)$$ is non-degenerated.
Let $e_1,\dots, e_d$ be a basis of $\Frac(\sO(X))$ over $\Frac \sO(\D^r)$.  
After multiplying by some element in $\sO(X)$, we may assume that $e_i\in \sO(X), i=1,\dots,d.$
Moreover, since $X$ is strict, we may assume that $\rho(e_i)=1, i=1,\dots, d.$
Since $X$ is irreducible and  the field extension $\pi^*: \Frac \sO(\D^r)\hookrightarrow \Frac(\sO(X))$ is separable,
$D:=\det((\Tr(e_ie_j)_{0\leq i,j\leq d}))\in \sO^{\circ}(\D^r)\setminus \{0\}.$

By \cite[3.2.8]{Berkovich1993}, there is a Zariski closed subset $Z'$ of $X$ such that $\pi$ is flat on $X\setminus Z'.$
By \cite[3.2.7]{Berkovich1993}, the map $\pi|_{X\setminus Z'}$ is open.
Set $Z:=\pi(Z')$. Since $\pi$ is finite and $X$ is irreducible, $Z$ is a Zariski closed subset of $Z$ and $\pi(X)\neq Z.$
There is $Q'\in \bk^{\circ}\{T_1,\dots,T_r\}$ such that $Q'|_Z=0$ and $Q'|_{\pi(X)}\neq 0.$
Then $Q:=\pi^*Q'\in \sO^{\circ}(X)\setminus \{0\}.$
By Lemma \ref{lemmultiqanna}, we only need to prove $QF\in \sO^{\circ}(X).$

\medskip

Define $W_Z^n:=\{x\in \D^r|\,\, |Q'(x)|\leq 1/2^n\}$, which is a Weierstrass subdomain of $\D^r$.
In particular $Z\subseteq (W_Z^n)^{\circ}.$

For every $x\in \D^r\setminus Z$, since $\pi^{-1}(x)$ is finite, there exists an open neighbourhood $V_x$ of $\pi^{-1}(x)$ such that every connected component of $V_x$ contains exactly one point $t\in \pi^{-1}(x).$
We denote such connected component by $V_x^t.$
For every $x\in \D^r, n\geq 0$, there exists a subset $J^n_x$ of $I^n$ such that 
$\pi^{-1}(x)\subseteq \cup_{i\in J_x^n}U^n_i.$ 
For $i\in I$, since $U^n_i$ is compact and $F^n_i\in \Frac(\sO(X))$, 
there is on open neighborhood $C^n_i$ 
of $U^n_i$ such that $F^n_i\in \sO(C^n_i)$ and $\rho(F^n_i|_{C^n_i}-F_n|_{C^n_i})\leq 1/2^{n-1}.$
There exists a map $\sigma_n:\pi^{-1}(x)\to I$ sending $t\in \pi^{-1}(x)$ to the a $j\in I$ containing $t.$
Since $\pi|_{X\setminus \pi^{-1}(Z)}$ is open, $(\cap_{t\in \pi^{-1}(x)} \pi(C^n_{\sigma^n(t)}))\cap (\cap_{t\in \pi^{-1}(x)} \pi(V_x^t))$ is a neighborhood of $x$.
Pick a rational neighborhood $W^n_x$ of $x$ contained in $(\cap_{t\in \pi^{-1}(x)} \pi(C^n_{\sigma^n(t)}\cap V_x^t))\cap (\D^r\setminus Z)$.
Then we have 
\begin{points}
\item $\pi^{-1}(W_x^n)\subseteq V_x;$
\item for every $t\in \pi^{-1}(x),$ $\rho((F^n_{\sigma_n(t)}-F_n)|_{W_x^n(t)})\leq 1/2^{n-1},$ where $W_x^n(t):=\pi^{-1}(W_x^n)\cap V_x^t$.
\end{points}
Then  we have $$\pi^{-1}(W_x^n)=\sqcup_{t\in \pi^{-1}(x)}W_x^n(t).$$
Then $W_x^n(t), t\in \pi^{-1}(x)$ are rational subdomains of $X$ and $$\sO(\pi^{-1}(W_x^n))=\prod_{t\in \pi^{-1}(x)}\sO(W_x^n(t)).$$
Then there exists $H^n_x\in \sO(\pi^{-1}(W_x^n))$ such that 
$H^n_x|_{W_x^n(t)}=QF^n_{\sigma_n(t)}.$
Then we have $$\rho(H^n_x-QF_n|_{\pi^{-1}(W_x^n)})\leq 1/2^{n-1}.$$
We note that $$\rho(QF_n|_{\pi^{-1}(W_Z^n)})\leq \rho(Q|_{\pi^{-1}(W_Z^n)})\leq 1/2^{n}.$$
Since $\D^r$ is compact, there is a finite subset $S^n\subset \D^r$ such that $X=(\cup_{x\in S^{n+1}}W_x^{n+1})\cup W_Z^n.$
After replacing $QF_n$ by $\{(\pi^{-1}(W_x^{n+1}),H^{n+1}_x), x\in S\}\cup \{(\pi^{-1}(W_Z^n),0)\}$, we may assume that for every $n\geq 0, i\in I^n$, $U^n_i$ takes form 
$U^n_i=\pi^{-1}(W^n_i)$, where $W_i^n$ is a rational subdomain of $\D^r$. Since we only need to prove $QF\in \sO^{\circ}(X),$ we may replace $F_n$ by $QF_n.$

For every $n\geq 0, i\in I^n, j=1,\dots, d,$ define $$H^{n,j}_{i}:=\Tr(e_jF^n_i)\in \Frac(\sO(\D^r))\cap \sO(W^n_i).$$
Define $H^{n,j}:=\{(W^n_i, H^{n,j}_{i}), i\in I^n\}.$ Then for every $m\geq n\geq 0$ we have 
$$\rho(H^{n,j})\leq \rho(F_n)\leq 1$$ and 
$$\rho(H^{n,j}-H^{m,j})\leq \rho(F_n-F_m)\leq 1/2^n.$$
Then we have $[H^{n,j}]\in \hat{\sO}^1(\D^r).$
By Lemma \ref{lemdnsectionglo}, there exists $B_n^j\in \sO^{\circ}(\D^r)$ such that $\rho(B_n^j-H^{n,j})\leq 1/2^n$.
So $B_n^j$ is a Cauchy sequence in $\sO^{\circ}(\D^r)$ and denote by $B^j$ its limit. 
We identify $\sO^{\circ}(\D^r)$ with its image in $\hat{\sO}^1(\D^r)$, we have $[H^{n,j}]=B^j.$

There exists a unique $G\in \Frac(\sO(X))$ such that $\Tr(e_jG)=B^j$ for all $j=1,\dots, r.$
\begin{lem}\label{lemregug}We have $G\times (\pi^*D)\in \sO^{\circ}(X).$
\end{lem}
\proof[Proof of Lemma \ref{lemregug}]
For every $x\in \D^r\setminus \{|D(x)|=0\}$, the morphism 
$\Phi_x: \sO(X)\hotimes_{\bk}\sH (x)\to {\sH}(x)^d$ sending $g$ to $(\Tr(e_i(x)g))_{i=1,\dots, s}$
is a homeomorphism.  On $\sO(Y)\hotimes_{\bk}\sH (x)$, we use the spectral norm $\|\cdot \|'_x$ and on ${\sH}(x)^d$ we use the usual product norm $\|\cdot \|_x$.
Then we have 
$$\|\Phi_x(g)\|_x'\leq \|g\|_x \text{ and } \|\Phi^{-1}(h)\|'_x\leq |D(x)|^{-1}\|h\|_x.$$

For every $x\in \D^r\setminus \{|D(x)|=0\}$,  $G$ is well defined on the fiber $\pi^{-1}(x)$ and
$$G(x)=\Phi^{-1}(B^1,\dots, B^d).$$
In particular, we have 
$$\max_{y\in \pi^{-1}(x)}|G(y)|=\|G(x)\|_x'\leq |D(x)|^{-1}.$$
It follows that for every $y\in X\setminus \{|D(x)|=0\}$ we have $|G\times (\pi^*D)(y)|\leq 1$.
Since $G\times (\pi^*D)\in \Frac(\sO(X))$ and $X$ is normal, we have $G\times (\pi^*D)\in \sO(X)$ and $\rho(G\times (\pi^*D))\leq 1.$
This concludes the proof.
\endproof

For every $x\in \D^r, n\geq 0, c\geq 1/2^n$, consider the function $\psi_{F_n}^c: \D^r\to \R$ defined by 
$$x\mapsto \max\{\max_{y\in \pi^{-1}(x)}\log|F^n_i(y)|, \log c\}=\max\{\log|F^n_i(x)|'_x, \log c\},$$ where $x\in W^n_i.$ We note that the value $\psi_{F_n}^i(x)$
does not depend on the choice of such $i.$ 
For every $x\in \D^r$, consider the character polynomial 
$$t^d+A^{n,d-1}_i(x)t^{n-1}+\dots+A^{n,0}_i(x)$$
of $F^n_i\in \sO(X)\hotimes_{\bk}\sH (x).$
Then for $x\in W^n_i$, we have that $$\psi_{F_n}^c(x)=\max\{\max_{j=1}^{d}(j^{-1}\log|A^{n,d-j}_i(x)|), \log(c)\}$$$$=\max_{j=1}^{d}(j^{-1}\max\{\log|A^{n,d-j}_i(x)|, j\log(c)\}\}.$$
Observe that $\{(W^n_i, A^{n_i,d-j}), i\in I^n\}\in \sP^1_{(1/2^n)^j}(\D^r).$
Let $\eta$ be the Gauss point of $\D^r.$
Apply Inequality \ref{eqmax} for the functions $x\mapsto \max\{\log|A^{n,d-j}_i(x)|, j\log(c)\}$, we have that 
\begin{equation}\label{equmaxpm}
\psi_{F_n}^c(x)\leq \psi_{F_n}^c(\eta).
\end{equation}
We note that our argument shows that this Inequality applies to every element in $\sP_{b}(X)$, not only for $F_n.$

\begin{lem}\label{lemtdconv}We have $$\rho(\pi^*(D)F_n-\pi^*(D)G)\leq  \max\{1, |D(\eta)|^{-1}\}1/2^n, n\geq 0.$$
In particular, we have 
$(\pi^*D) F=[\pi^*D\times F_n]=G\times (\pi^*D)\in \sO^{\circ}(X).$
\end{lem}
\proof[Proof of Lemma \ref{lemtdconv}]
For every $x\in \D^r\setminus \{|D(x)|=0\}, n\geq 0, x\in W_i^n$, we have 
$$F^n_i(x)-G(x)=\Phi^{-1}((H_i^{n,1}-B^1,\dots,H_i^{n,d}-B^d)).$$
In particular, we have 
\begin{equation}\label{equpbufg}
\|F^n_i(x)-G(x)\|_x'\leq |D(x)|^{-1}/2^n.
\end{equation}
Observe that $|D(\eta)|>0$. So we have 
$\|F^n_i(\eta)-G(\eta)\|_x'\leq |D(\eta)|^{-1}/2^n.$
Then we have 
$$\|\pi^*(D)F^n_i(\eta)-\pi^*(D)G(\eta)\|_x'\leq |D(\eta)|^{-1}1/2^n.$$
Since $F_n\in \sP^1_{1/2^n}(X)$ and $\pi^*D,G\in \sO^{\circ}(X)$, 
$\pi^*(D)F_n-\pi^*(D)G\in  \sP^1_{1/2^n}(X).$
By Inequality \ref{equmaxpm}, we get 
$$\rho(\pi^*(D)F_n-\pi^*(D)G)\leq  \max\{1, |D(\eta)|^{-1}\}1/2^n,$$
which concludes the proof.
\endproof

Then we conclude the proof by Lemma \ref{lemmultiqanna}.
%

\subsection{Rigid space version of Theorem \ref{thmisolim}}\label{subsectionanotherproof}
In this section, we give a simplified proof for the rigid space version of Theorem \ref{thmisolim} via \cite[Theorem 2]{Bartenwerfer1978}.  This proof is suggested by the referee of the first version of our paper. 

By the discussions in the beginning of Section \ref{sectionproofthmmain} and in Section \ref{subsectionredtostrict}, 
we may assume that $X$ is a strict affinoid space and we only need to show that the morphism $\alpha_X: \sO^{\circ}(X)\to \lim\limits_{\leftarrow}\sO^{\circ}/\omega^n(X)$ is surjective.
Not the we only consider the strict covers in this section.

By rigid space version of Proposition \ref{prodxlimo}, we identify $\lim\limits_{\leftarrow}\sO^{\circ}/\omega^n(X)$ with $\hat{\sO}^1(X)$ and identify $\sO^{\circ}(X)$ as its image in $\hat{\sO}^1(X).$

Let $F=[F_n]$ be an element in $\hat{\sO}^1(X)$. We may assume that $F_n\in \sP^1(X).$
Write $F_n=\{(U_i^n, F^n_i), i\in I^n\}.$ Note that here we ask the $U_i$ are strict.
By the argument in Section \ref{subsectionreducefinite}, we may assume that $I^n$ are finite and $U_i^n$ are rational subdomains in $X$.
We may assume that for every $m\geq n\geq 0,$
$\rho(F_m-F_n)\leq 1/2^n.$

\medskip

We first treat the case where $X$ is smooth.  
As in \cite{Bartenwerfer1978}, for every strict affinoid cover $\sU$ of $X$ and $r>0$, let $C^{\d}_r(\sU)$ the C\v{e}ch complex  for $\sO_X^{<r}$ where $\sO_X^{<r}$ is the sheaf of analytic functions of spectral norm $<r$.  Let $$H_r^{\d}(X):=\lim_{\rightarrow}H_r^{\d}(\sU)$$ where the limit is taken over the system of all strict affinoid covers.
By \cite[Theorem 2]{Bartenwerfer1978}, there is $c\in \bk$ with $0<|c|\leq 1$ such that 
\begin{equation}\label{equrigidc}c\cdot H_r^1(X)=0.
\end{equation}
For every $n\geq 0$, $F_n$ defines a cocycle $F^n_{i,j}:=\{(F^n_i-F^n_j)|_{U_i^n\cap U_j^n}\}.$
By (\ref{equrigidc}), there is a cover $\sU'=\{U_i', i\in I'\}$ refines $\{U_i^n, i\in I^n\}$ via $\sigma:I'\to I$ and $G_i\in \sO^{<2^{n-1}}(U_i')$ such that 
$G_i-G_j=cF^n_{\sigma(i),\sigma(j))}$ on $U'_i\cap U'_j.$ Then there is $G_n\in \sO(X)$ such that $G_n|_{U'_i}=F^n_{\sigma(i)}-c^{-1}G_i.$
Then we have $\rho(G_n-F_n)\leq |c^{-1}|/2^{n-1}.$ It follows that $[F_n]=\lim_{n\to \infty}G_n\in \sO(X),$ which concludes the proof.

We now prove the general case. Pick $G\in \sO(X)$ such that $Y:=X\setminus \{G=0\}$ is smooth. 
The proof of the smooth implies that $F|_Y$ is contained in $\hat{\sO}^1(Y).$
Since $X$ is normal, we conclude the proof by Riemann extension theorem \cite[3.3.14]{Berkovich1990}.

\section{Analytic and continuous functions on subsets}
Let $X$ be a separated and reduced $\bk$-analytic space.
For every subsets $Y\subseteq Z\subseteq X$, we denote by $Y^{\circ,Z}$ to be the interior of $Y$ as a subset of $Z$. It is clear that $Y^{\circ}\subseteq Y^{\circ,Z}$.

Let $S$ be a subset of $X$. 
We say that a collection of analytic subdomains $U_i,i\in I$ is a \emph{$G$-covering of $S$} if for every point $x\in S$, there is an analytic subdomain $U$ and an open subset $W$ such that $x\in W\cap S\subseteq U$ and $U\cap U_i, i\in I$ is a $G$-covering of $U.$

\begin{rem}Open coverings of $S$ and finite coverings of $S$ by analytic subdomains are $G$-covering of $S$.
\end{rem}

\subsection{Functions on $S$}
Define $\Fun(S):=\prod_{x\in S}\sH(x)$ and we call an element $f=(f(x))_{x\in S}\in \Fun(S)$ a \emph{function} on $S$.
Define $\rho_S: \Fun(S)\to [0,\infty]$ by $\rho_S(f):=\sup_{x\in S}|f(x)|.$ Sometimes we write $\rho$ for $\rho_S$ for the simplicity.
It is clear that $\Fun(S)$ is $\bk$-algebra. For $f,g\in \Fun(S)$, we have 
\begin{points}
\item $f=g$ if and only if $\rho(f-g)=0;$
\item $\rho(f+g)\leq \max\{\rho(f),\rho(g)\};$
\item $\rho(fg)\leq \rho(f)\rho(g).$
\end{points}
The space $\Fun(S)$ is complete in the following sense:
For $f_n\in \Fun(S), n\geq 0$, if $\lim\limits_{n\to \infty}\rho(f_n)=0$, let $\sum_{n\geq 0}f_n$ the function such that for $x\in S$,  $(\sum_{n\geq 0}f_n)(x):=\sum_{n\geq 0}f_n(x)\in \sH(x).$
This is the unique function $g\in \Fun(S)$ such that $\rho(g-\sum_{i=0}^nf_n)\to 0$ as $n\to \infty.$

For $S\subseteq S'$, we have a natural restriction morphism $f=(f(x))_{x\in S'}\mapsto f|_S:=(f(x))_{x\in S}.$ We have $\rho(f|_S)\leq \rho(f).$

Let $U$ be an analytic subdomain of $X$ and $r\geq 0$. An elements in $F=\{(U_i,F_i), i\in I\}\in \sP(U)$ can be viewed as a multi-value function on $U.$ 
For $x\in U$ and $f(x)\in \sH(x)$, we denote by $|F(x)-f(x)|:=\sup_{i\in I}|F_i(x)-f(x)|$.
For a subset $S\subseteq U$ and $f\in \Fun(S)$, 
define $\rho(f-F|_S):=\sup_{ x\in S}|F(x)-f_i(x)|=\sup_{i\in I, x\in S}|F_i(x)-f_i(x)|.$

\begin{rem}\label{remtextone}For $x\in U, f(x)\in \sH(x), s\geq r\geq 0$ and $F=\{(U_i,F_i), i\in I\}\in \sP_r(U)$, $|F(x)-f(x)|\leq s$ if and only if there is $i\in I$ such that $x\in U_i$ and $|F_i(x)-f(x)|\leq s.$
\end{rem}

For every $x\in X$ and $f(x)\in \sH(x)$. For $r\geq   0$, an \emph{$r$-extension of $f(x)$} is a pair $(U, F)$ where $U$ is an analytic neighborhood of $x$  and $F\in \sP_r(U)$ such that $|F(x)-f(x)|\leq  r.$
We say that $(U, F)$ is \emph{finite affinoid}, if $F=\{(U_i, F_i), i\in I\}\in \sP(U)$, where $I$ is a finite, and $U_i, i\in I$ are affinoid.
\begin{lem}\label{lemexiuniqrex}
Let $U$ be an analytic neighborhood of a point $x\in X$.
Let $f(x)\in\sH(x)$ and $r>0$. 
For $F=\{(U_i, F_i), i\in I\}\in \sP(U)$, if $I$ is a finite, $U_i, i\in I$ are compact and $|F(x)-f(x)|<r$, then there is an open neighborhood $U'$ of $x$ contained in $U$ such that $F|_U'$ is an $r$-extension of $x$.

In particular, for every $x\in X$, $f(x)\in\sH(x)$ and $r>0$,  
there exists a finite affinoid $r$-extension $F=\{(U_i, F_i), i\in I\}$ of $x$.

Moreover, for two $r$-extensions $(V_1,F_1)$ and $(V_2,F_2)$ of $f(x)$ and $s>r$, there is an open neighborhood $W$ of $x$ such that $W\subseteq V_1\cap V_2$ such that $\rho(F_1|_W-F_2|_W)< s.$
\end{lem}
\proof[Proof of Lemma \ref{lemexiuniqrex}]
After restricting $F$ on a open neighborhood of $x$, we may assume that $x\in (\cap_{i\in I}U_i)\cap (\cup_{i\in I}U_i)^{\circ}.$
For every $i,j=1,\dots,m$, we have $F_i|_{U_i\cap U_j}-F_j|_{U_i\cap U_j}\in \sO(U_i\cap U_j).$  Since $|(F_i|_{U_i\cap U_j}-F_j|_{U_i\cap U_j})(x)|< r,$
there is an open neighborhood $U'$ of $x$ such that for every $i,j=1,\dots,m$ and $y\in U_i\cap U_j\cap U'$, $|(g_i|_{U_i\cap U_j}-g_j|_{U_i\cap U_j})(y)|<r.$
Then $F|_{U'}\in \sP_r((\cup_{i=1}^mU_i')\cap U)$ is an $r$-extension of $f(x).$

Let $x$ be a point in $X$.
There are finite affinoid subdomains $U_1,\dots, U_m$ of $X$ such that $x\in (\cap_{i=1}^mU_i)\cap (\cup_{i=1}^mU_i)^{\circ}.$
For every $i=1,\dots,m$, there is an affinoid subdomain $U_i'$ of $U_i$ such that $x\in U_i'^{\circ,U_i}$ and $g_i\in U_i'$ such that $|g_i(x)-f(x)|\leq r/2.$
We note that $U:=\cup_{i=1}^mU_i'$ is an analytic neighborhood of $x.$ 
Apply the previous paragraph for $F:=\{(U_i', g_i), i=1,\dots,m\}\in \sP(U)$, is an open neighborhood $U'$ of $x$ contained in $U$ such that $F|_U'$ is an $r$-extension of $x$.
There are affinoid subdomains $V_1,\dots, V_l$ of $U'$ such that $x\in (\cup_{i=1}^l V_i)^{\circ}$. Then $F':=\vee_{i=1}^lF|_{V_i}$ is an $r$-extension of $f(x)$ as we need.


Let $(V_1,F_1)$ and $(V_2,F_2)$ be two $r$-extensions of $f(x).$ We have $F_1-F_2\in \sP_r(V_1\cap V_2)$ and $|F_1(x)-F_2(x)|\leq r<s.$
By Proposition \ref{procontinuityphif}, the function $\Phi:V_1\cap V_2\to [r, \infty)$ sending $y$ to $\max\{|F_1(y)-F_2(y)|, r\}$ is a continuous function. 
Since $\Phi(x)=r<s$,  there is an open neighborhood $W$ of $x$ such that $W\subseteq V_1\cap V_2$ such that
 $\rho(F_1|_W-F_2|_W)< s.$
\endproof

\subsection{$r$-Analytic functions on $S$}
Let $U$ be an analytic subdomain of $X$ containing $S$, and $f\in \sO(U)$, denote by  $f(x)$ the image of $f$ in $\sH(x)$, then the morphism $f\mapsto (f(x))_{x\in U}\in \Fun(U)$ gives a natural embedding 
$\sO(U)\subseteq \Fun(U).$

For $r\in [0,\infty)$, an element $f\in \Fun(S)$ is said to be \emph{$r$-analytic} if there is a $G$-covering $U_i, i\in I$ of $S$, $f_i\in \sO(U_i)$ such that 
$\rho(f|_{S\cap U_i}-f_i|_{S\cap U_i})\leq r.$ Define $\Fun^{r-\an}(S)$ to be the set of $r$-analytic functions on $S$.  
For $f,g\in \Fun(S)$, if $\rho(f-g)\leq r$ then $f\in \Fun^{r-\an}(S)$ if and only if $g\in \Fun^{r-\an}(S).$

Define $\Fun^{\an}(S):=\Fun^{0-\an}(S)$. The elements in $\Fun^{\an}(S)$ are called the \emph{analytic functions} on $S$.
It is clear that for $S\subseteq S'$ and $f\in \Fun^{\an}(S')$, $f|_S\in \Fun^{\an}(S).$

\begin{rem}
When $S$ is an analytic subdomain,  then $\Fun^{\an}(S)=\sO(S).$
More generally, when $X$ is a Zariski closed subset of some analytic subdomain of $X$, we have $\Fun^{\an}(S)=\sO(S).$ 
\end{rem}


The following property shows that 
to be $r$-analytic is a local property.
\begin{pro}\label{proanalocal}For $f\in \Fun(S)$, the following statements are equivalent:
\begin{points}
\item $f$ is $r$-analytic on $S$; 
\item for every analytic subdomain $V\subseteq X$, $f|_{S\cap V}$ is $r$-analytic on $S\cap V;$
\item for every $x\in S$, there is an open neighborhood $V$ of $x$, such that $f$ is $r$-analytic on $S\cap V;$
\item for every $x\in S$, there is an open neighborhood $V$ of $x$, an analytic subdomain $W$ in $V$ containing $S\cap V$ and $F\in \sP(W)$ such that 
$\rho(F|_S-f|_S)\leq r;$
\item for every $s>r, x\in S$, there is an open neighborhood $V$ of $x$, an analytic subdomain $W$ in $V$ containing $S\cap V$ and $F\in \sP_s(W)$ such that 
$\rho(F|_S-f|_S)\leq r.$
\end{points}
\end{pro}

\proof[Proof of Proposition \ref{proanalocal}]
It is clear that (i) implies (ii), (ii) implies (iii) and (v) implies (i).

We first prove that (iii) implies (iv). 
For $x\in S$, there is an open neighborhood $U$ of $x$, such that $f|_{S\cap U}$ is analytic on $S\cap U.$
There is a $G$-covering $V_i, i\in I$ of $S$, $f_i\in \sO(V_i)$ such that 
$\rho(f|_{S\cap V_i}-f_i|_{S\cap V_i})\leq r.$
Then there is an open neighborhood $V$ of $x$, an analytic subdomain $W$ in $V$ containing $S\cap V$ such that $W\cap V_i,i\in I$ is a $G$-covering of $S\cap W.$
Set $F:=\{(W\cap V_i, f_i|_{W\cap V_i}), i\in I\}\in \sP(W)$. Then we have $\rho(F|_S-f|_S)\leq r.$

Now we only need to prove that (iv) implies (v).  Assume that there is an open neighborhood $V$ of $x$, an analytic subdomain $W$ in $V$ containing $S\cap V$ and $F=\{(U_i,F_i), i\in I\}\in \sP(W)$ such that 
$\rho(F|_S-f|_S)\leq r.$
There is a finite subset $J$ of $I$ such that $x\in (\cap_{i\in J}U_i)\cap (\cup_{i\in J}U_i)^{\circ, W}.$
For every $i\in J$,  $|f(x)-F_i(x)|\leq r<s$. So for every $i,j\in I$, $|F_i(x)-F_j(x)|\leq r<s$.
By Proposition \ref{procontinuityphif}, there is an analytic subset $V$ in $W$ such that  $x\in V^{\circ, W}$ and for $i,j\in J$,
$\rho_{V_i\cap V_j}(F_i-F_j)<s.$ After replacing $W$ by $V$ and $F$ by $\{(V_i\cap V, F_i|_{V_i\cap V}), i\in J\}$, we conclude the proof.
\endproof

By Proposition \ref{procontinuityphif}, for every $r\geq 0$, $f\in \Fun^{r}(S),$ the map $\phi_f^r: S\to [-\infty,\infty)$ sending 
$x$ to $\max\{\log |f(x)|, \log|r|\}$ is continuous.

\subsection{Continuous functions}Define $\Fun^{\Con}(S):=\cap_{r>0}\Fun^{r-\an}(S)$. The elements in $\Fun^{\Con}(S)$ are called the \emph{continuous functions} on $S$.
It is clear that $\Fun^{\an}(S)\subseteq \Fun^{\Con}(S)$ and  $\Fun^{\Con}(S)$ is a $\bk$-subalgebra of $\Fun(S)$. Moreover, the following property shows that $\Fun^{\Con}(S)$ is complete.
\begin{pro}\label{proconfuncomp}
For $f_n\in \Fun^{\Con}(S), n\geq 0$ with $\lim\limits_{n\to \infty}\rho(f_n)= 0$, we have $\sum_{n\geq 0}f_n\in \Fun^{\Con}(S).$
\end{pro}
\proof[Proof of Proposition \ref{proconfuncomp}]
Set $g_m:=\sum_{n=0}^mf_n$ and $g:=\sum_{n\geq 0}f_n$.
For every $\epsilon>0$, there is $N\geq 0$ such that for $n\geq N$, $\rho(f_n)\leq \epsilon.$
Then we have $\rho(g-g_N)\leq \epsilon.$
Since $g_N\in \Fun^{\Con}(S)\subseteq \Fun^{\epsilon-\an}(S),$ we have $g\in \Fun^{\epsilon-\an}(S).$
Then we get $g\in \cap_{\epsilon>0}\Fun^{\epsilon-\an}(S)=\Fun^{\Con}(S).$
\endproof

When $S$ is an analytic subdomain of $X$, we have $\Fun^{\Con}(S)=\hat{\sO}(S)$.
More generally, when $X$ is a Zariski closed subset of some analytic subdomain of $X$, we have $\Fun^{\Con}(S)=\hat{\sO}(S).$ 
Remark \ref{remnormalness} shows that, in general,  $\hat{\sO}(S)$ is not dense in $\Fun^{\Con}(S).$

By Theorem \ref{thmisolim}, when $S$ is an analytic subdomain of $X$, if $\bk$ is perfect and $S$ is normal, then $\Fun^{\Con}(S)=\sO(S)=\Fun^{\an}(S).$ 
When $S=\{x\}$ is a single point in $X$, we have $\Fun^{\Con}(\{x\})=\sH(x)$ which is the completion of $\Fun^{\an}(\{x\})$.
The following example shows that in general $\Fun^{\an}(\{x\})\neq \Fun^{\Con}(\{x\}).$
\begin{exe}
Assume that $\bk=\C((t)).$ Let $\xi$ be the Gauss point of $\D=\sM(\bk\{T\}).$ Set $f:=\sum_{i\geq 0} \frac{t^i}{T+i}$.
Then we have $f\in \sH(\xi)\setminus \Fun^{\an}(\{\xi\}).$
\end{exe}

\medskip

For $x\in S$, and $f\in \Fun(S)$, we say that $f$ is \emph{continuous at $x$} if for every $r>0$, there is an $r$-extension $(U,F)$ of $x$
such that $\rho(F|_{S\cap U}-f|_{S\cap U})\leq r$.

\begin{pro}\label{procricon}
For $f\in \Fun(S)$, the following statements are equivalent:
\begin{points}
\item $f$ is continuous;
\item $f$ is continuous at every point $x\in S.$
\end{points}
\end{pro}
\proof[Proof of Proposition \ref{procricon}]
By Proposition \ref{proanalocal}, it is clear that (ii) implies (i).

We only need to show that (i) implies (ii).
Since $f$ is continuous, by Proposition \ref{proanalocal}, for every $r>0,$
there is an open neighborhood $V$ of $x$, an analytic subdomain $W$ in $V$ containing $S\cap V$ and $F=\{(V_i, F_i), i\in I\}\in \sP_{r/2}(W)$ such that 
$\rho(F|_S-f|_S)\leq r/3.$ After replacing $V_i, i\in I$ by refinement by an affinoid covering, we may assume that $V_i, i\in I$ are affinoid. 
There is a finite subset $J\subseteq I$ such that $x\in W':=\cup_{i\in J}V_i.$
Set $F':=\{(V_i, F_i), i\in J\}.$ By Lemma \ref{lemexiuniqrex}, we may pick a finite affinoid $r/2$-extension $(U,G)$ of $f(x).$
Then we have $H:=G\vee F'|_{W'\cap U}\in \sP(U)$ and $|H(x)-f(x)|\leq r/2.$ 
By Lemma \ref{lemexiuniqrex}, there is an open neighborhood $U'$ of $x$ such that $U'\subseteq U$ and 
 $H|_{U'}\in \sP_{r}(U).$ Then $(U', H|_{U'})$ is an $r$-extension of $f(x).$ 
 By  Remark \ref{remtextone}, we have $\rho(H|_{U'\cap S}-f|_{U'\cap S})\leq r.$ So $f$ is continuous at $x$.
\endproof

For $f\in \Fun(S)$, we define a map $\Phi_f: S\to \A^{1,\an}=\sM(\bk[T])$ sending $x\in S$ to the point in $\A^{1,\an}$ defined by the morphism $\bk[T]\to \sH(x)$ sending $T\to f(x).$

\begin{pro}\label{proconinducon}For $f\in \Fun^{\Con}(S)$, the induced map $\Phi_f$ is continuous.
\end{pro}

\proof[Proof of Proposition \ref{proconinducon}]
For every $g\in \bk[T]$, we have $g\circ f\in \Fun^{\Con}(S).$ Since for every $r>0$, $g\circ f\in \Fun^{r-\an}(S)$, the map 
$\phi_{g\circ f}^r: S\to \R$ sending $x$ to $\max\{\log|g\circ f(x)|, \log r\}$ is continuous. 
Let $r$ tend to $0$, the function $\phi_{g\circ f}: x\mapsto \log|g\circ f(x)|=\log|g\circ (\Phi_f)|$ is continuous for every $g\in \bk[T]$. This shows that $\Phi_f$ is continuous.
\endproof

The inverse of Proposition \ref{proconinducon} does not hold.
\begin{exe}Assume that $\bk=\C((t))$. Let $X:=\D=\sM(\bk\{T\})$ and $S:=\Z\cup \{\xi\}$ where $\xi$ is the Gauss point of $\D.$
It is clear the $S$ is a closed subset of $X$.
Define $f\in \Fun(S)$ by $f(i):=i\in \Z\subseteq \bk$ and $f(\xi):=T^2.$
Then $\Phi_f$ is the identity map on $S$, so it is continuous.
But $\Phi_f$ is not continuous, because $\Phi_{f-T|_S}$ is not continuous. 
\end{exe}

\rem
More generally, give finitely many continuous functions $f_1,\dots, f_r\in \Fun(S)$, they naturally  induce a continuous map 
$$\Phi_{f_1,\dots, f_r}: S\to \A^{r,\an}=\sM(\bk[T_1,\dots, T_r]),$$
defined by the morphism $\bk[T_1,\dots,T_r]\to \sH(x)$ sending $T_i\to f_i(x).$
Since the underlining topological space of $\A^{r,\an}$ is not the product of the underlining topological space of $\A^{1,\an}$,
usually we do not has a natural way to get a ``product" from $r$ continuous maps from $S$ to $\A^{1,\an}$. This suggest that our definition of continuous functions on $S$ should be the natural one.  Since every affinoid space is a subset of some $\A^{r,\an}$, it is not hard to define continuous maps between subsets of analytic spaces via our notion of continuous functions. 
\endrem

For every $x\in X(\bk)$, we have $\sH(x)=\bk$.
So when $S$ is a closed subset of $X(\bk)$, a function $f\in \Fun(S)$ can be viewed a map $f: S\to \bk.$
The follows proposition shows that when $S\subseteq X(\bk)$ our definition of continuity coincides the usual one.
\begin{pro}\label{proclosedpointcon}Let $S$ be a subset of $X(\bk)$. A function $f\in \Fun(S)$ is continuous if and only if the induced map $f: S\to \bk$ is continuous.
\end{pro}
\proof[Proof of Proposition \ref{proclosedpointcon}]We identify $\bk$ with $\A^{1,\an}(\bk)$. By Proposition \ref{proconinducon}, if $f$ is continuous, the map $f: S\to \bk\subseteq \A^{1,\an}(\bk)$ is continuous.

Now we assume that the map $f: S\to \bk$ is continuous. 
Let $x$ be a point in $S.$ For every $r>0$, there is an open neighborhood $U$ of $x$ such that for every $y\in S\cap U$, $|f(y)-f(x)|<r.$
We note that $f(x)\in \bk$ can be viewed as a constant function on  $X$.
Then $\{(U, f(x)|_U)\}$ is an $r$-extension of $f(x)$. Then $f$ is continuous at $x$.
We conclude the proof by Proposition \ref{procricon}.
\endproof

In the end, we give an example of a continuous function on some closed set in $\D(\bk)$ which is everywhere non-analytic.
\begin{exe}
Assume that $\bk=\Q_p$ where $p$ is a prime. Set $S:=\Z_p:= \D(\Q_p)$ which is a closed subset of $X:=\D=\sM(\bk\{T\}).$ We may view $S$ as a Cantor set in $\D.$

Every $x\in \Z_p$ can be uniquely written as $\sum_{i\geq 0}x_ip^i$ where $x_i\in \{0,\dots, p-1\}.$
Let $f: \Z_p\to \bk$ be the function $\sum_{i\geq 0}x_ip^i\mapsto \sum_{i\geq 0}x_ip^{2i}.$
This function is non-constant on every open subset of $\Z_p.$
For $x,y\in \Z_p$, we have $|f(x)-f(y)|=|x-y|^2.$
By Proposition \ref{proclosedpointcon}, $f$ is continuous. 

We claim that $f$ is everywhere non-analytic i.e. for every nonempty open subset $U$ of $\Z_p$
$f|_U\not\in \Fun^{\an}(U).$
Otherwise, there is $x\in \Z_p$ and an affinoid subset $U$ of $\D$ containing $x$ and a function $g\in \sO(U)$ such that $g|_{U\cap \Z_p}=f.$ 
Since $x\in U(\bk)$, we have $x\in U^{\circ}.$ There is $r>0$ such that the disk $\D(x,r):=\{y|\,\, |(T-x)(y)|\leq r\}$ is contained in $U$.
We may assume that $U=\D(x,r)$ and $g$ takes form $\sum_{i\geq 0} a_i(T-x)^i$ where $a_i\in \Q_p$ and $|a_i|r^i\to 0$ as $i\to \infty.$
Since for every $y,z\in \Z_p\cap U$, $|g(y)-g(z)|=|f(y)-f(z)|=|y-z|^2$, $g'|_{U\cap \Z_p}=0.$ Since $U\cap \Z_p$ is Zariski dense in $U$,
we get $g'=0.$ So $g$ is constant in on $U$. It implies that $f$ is constant on $U\cap S$, which is a contradiction. 
\end{exe}

\bibliography{dd}
\end{document}